\documentclass[10pt]{amsart}
\usepackage{amsmath,amscd,amssymb}

\newtheorem{theorem}{Theorem}[section]

\newtheorem{proposition}[theorem]{Proposition}

\newtheorem{lemma}[theorem]{Lemma}
\theoremstyle{definition}    
\newtheorem{definition}[theorem]{Definition}
\theoremstyle{remark}

\newtheorem{remark}[theorem]{Remark}

\newcommand{\pt}{\operatorname{pt}}

\newcommand\A{\mathcal{A}}

\newcommand\G{\mathcal{G}}
\newcommand{\K}{\mathbb{K}}

\renewcommand{\L}{\mathcal{L}}
\renewcommand{\O}{\mathcal{O}}

\newcommand{\Co}{\mathcal{C}}

\newcommand{\ca}{\mathcal}

\newcommand{\U}{\on{U}}
\newcommand{\E}{\ca{E}}
\renewcommand{\S}{\mathcal{S}}

\newcommand{\F}{\mathcal{F}}

\newcommand{\R}{\mathbb{R}}
\newcommand{\C}{\mathbb{C}}
\newcommand{\SU}{\on{SU}}

\newcommand{\Z}{\mathbb{Z}}

\newcommand\lie[1]{\mathfrak{#1}}

\newcommand{\g}{\lie{g}}

\newcommand{\z}{\lie{z}}

\renewcommand{\t}{\lie{t}}

\newcommand{\on}{\operatorname}

\newcommand{\Ad}{ \on{Ad} }

\newcommand{\ad}{\on{ad}}

\renewcommand{\ker}{ \on{ker}}

\newcommand{\Spin}{ \on{Spin}}

\newcommand{\SO}{ \on{SO}}

\newcommand{\vol}{  \on{vol}}

\newcommand{\tpi}{{2\pi i}}

\newcommand\qu{/\kern-.7ex/} 

\renewcommand{\d}{{\mbox{d}}}
\newcommand{\ol}{\overline}

\newcommand\Phinv{\Phi^{-1}}

\newcommand\Sig{\Sigma}

\newcommand\Om{\Omega}
\newcommand\om{\omega}

\newcommand{\f}{\frac}

\newcommand{\p}{\partial}

\newcommand\hh{{\f{1}{2}}}
\newcommand{\ti}{\tilde}
\newcommand{\eeq}{\end{eqnarray*}}
\newcommand{\beq}{\begin{eqnarray*}}

\renewcommand{\H}{\ca{H}}

\newcommand{\Cl}{{\C\on{l}}}
\newcommand{\pr}{\on{pr}}
\newcommand{\wh}{\widehat}

\newcommand{\mf}{\mathfrak}
\newcommand{\su}{\mf{su}}
\newcommand\dirac{/\kern-1.2ex\partial} 
\newcommand{\reg}{{\on{reg}}}
\begin{document}

\title{Quantization of q-Hamiltonian SU(2)-spaces}
\author{E. Meinrenken}
\address{University of Toronto, Department of Mathematics,
40 St George Street, Toronto, Ontario M4S2E4, Canada }
\email{mein@math.toronto.edu}

\date{\today}
\begin{abstract}
  We will explain how to define the quantization of q-Hamiltonian
  $\SU(2)$-spaces as push-forwards in twisted equivariant
  $K$-homology, and prove the `quantization commutes with reduction'
  theorem for this setting. As applications, we show how the Verlinde
  formulas for flat $\SU(2)$ or $\SO(3)$-bundles are obtained via
  localization in twisted $K$-homology.
\end{abstract}

\maketitle

\begin{quote}
{\it \small Dedicated to Hans Duistermaat on the occasion of his 65th birthday.}
\end{quote}
\vskip1cm

\section{Introduction}
The theory of q-Hamiltonian $G$-spaces was introduced ten years
ago in the paper \emph{Lie group-valued moment maps} \cite{al:mom}.
The motivation was to treat Hamiltonian loop group actions
with proper moment maps in a purely finite-dimensional framework,
obtaining for instance a finite-dimensional construction of the moduli
space of flat $G$-bundles over a surface. Many of the standard
constructions for ordinary Hamiltonian group actions on symplectic
manifolds carried over to the new setting, but often
with non-trivial `twists'. For example, all q-Hamiltonian $G$-spaces
$M$ carry a natural volume form \cite{al:du}, which may be viewed informally as a
push-forward of the (ill-defined) Liouville form on the associated
infinite-dimensional loop group space. This volume form admits an
equivariant extension (but for a non-standard equivariant cohomology
theory) \cite{al:gr}, and the total volume may be computed by
localization techniques, just as in the usual Duistermaat-Heckman
theory \cite{du:on}.

One problem that had remained open until recently is how to define a
`quantization' of q-Hamiltonian spaces. In contrast to the Hamiltonian
theory, the 2-form on a q-Hamiltonian space is usually degenerate.
Hence, there is no obvious notion of a compatible almost complex
structure, and the usual quantization as the equivariant index of a
$\Spin_c$-Dirac operator \cite{du:he} is no longer possible. In a
forthcoming paper \cite{al:for}, rather than trying to construct such
an operator, we define the quantization more abstractly as the
push-forward of a $K$-homology fundamental class $[M]$.  This
fundamental class is canonically defined as an element in
\emph{twisted} equivariant $K$-homology of $M$. Our construction
defines a push-forward of this element to the twisted equivariant
$K$-homology of a Lie group.  The Freed-Hopkins-Teleman theorem
\cite{fr:lo1,fr:tw} identifies the latter with the fusion ring
$R_k(G)$ (Verlinde algebra), at an appropriate level $k$. We take the
resulting element $\ca{Q}(M)\in R_k(G)$ to be the `quantization' of
our q-Hamiltonian space. As in the usual Hamiltonian theory
\cite{gu:ge,gu:re,me:sym}, the quantization procedure satisfies a
`quantization commutes with reduction' principle.

In the present paper, we will preview this quantization of
q-Hamiltonian $G$-spaces for the simplest of simple compact Lie
groups $G=\SU(2)$. Much of the general theory simplifies in this
special case -- for example, there is a fairly simple proof of the
q-Hamiltonian `quantization commutes with reduction' theorem. As an
application, we explain, following \cite{al:fi}, how the
$\SU(2)$-Verlinde formulas are obtained in our theory. In the last
Section, we will show how to derive Verlinde-type formulas for moduli
spaces of flat $\SO(3)$-bundles. The paper will be largely
self-contained, except for certain details that are
better handled with the techniques from  \cite{al:for}.\\[0.7cm]

\noindent{\bf Notation.}
We fix the following notations and conventions for the Lie group
$\SU(2)$.  The group unit will be denoted $e$, and the non-trivial central
element $c=\on{diag}(-1,-1)$. We define an open cover by contractible subsets
\begin{equation}\label{eq:cover}
 \SU(2)_+=\SU(2)\backslash \{c\},\ \
\SU(2)_-=\SU(2)\backslash\{e\}\end{equation}
with intersection the set $\SU(2)_\reg$ of regular elements. We take
the maximal torus $T$ to consist of the diagonal matrices, isomorphic 
to $\U(1)$ by the homomorphism
\[ j\colon \U(1)\to T,\ z\mapsto \on{diag}(z,z^{-1}).\]
The Weyl group $W=\Z_2$ acts on $T$ by permutation of the diagonal
entries, or equivalently on $\U(1)$ by $z\mapsto z^{-1}$. We let
$\Lambda\subset \t$ be the integral lattice (kernel of $\exp|_\t$) and
$\Lambda^*\subset \t^*$ its dual, the (real) weight lattice. For any
$\mu\in\Lambda^*$ we denote by $t\mapsto t^\mu$ the corresponding
homomorphism $T\to \U(1)$; the resulting 1-dimensional representation
of $T$ is denoted $\C_\mu$.  The weight lattice is generated by the
element $\rho\in\Lambda^*$ such that $\C_\rho$ is the defining
representation of $\U(1)$.  The corresponding positive root is
$\alpha=2\rho$. We will identify $\su(2)^*\cong\su(2)$ using the
\emph{basic inner product}
\[\xi\cdot\xi'=\f{1}{4\pi^2}\on{tr}(\xi^\dagger \xi'), \ \ \
\xi,\xi'\in\mf{su}(2).\] 
Similarly we identify $\t\cong\t^*$ using the induced inner product.
Under this identification, $\Lambda=2\Lambda^*$, with generators 
$\alpha=\tpi \on{diag}(1,-1)$ and $\rho=i\pi \on{diag}(1,-1)$. 

For any subset $A\subset \t$, we denote $T_A=\exp(A)=\{\exp\xi|\ 
\xi\in A\}$. Any conjugacy class in $\SU(2)$ passes through a unique
point in $T_{[0,\rho]}$, so that $[0,\rho]$ labels the conjugacy
classes.  We will frequently use the equivariant diffeomorphism,
\begin{equation}\label{eq:useful}
 T_{(0,\rho)}\times \SU(2)/T\to \SU(2)_\reg,\ (t,gT)\mapsto
\Ad_g(t).\end{equation}

\section{The fusion ring $R_k(\SU(2))$}
In this Section, we review three simple descriptions of the level $k$
fusion ring (Verlinde algebra) $R_k(G)$ for the case $G=\SU(2)$. The
fusion ring may be identified with the set of irreducible projective
representations of the loop group $L\SU(2)$ at level $k$ \cite{pr:lo}, but we will
not need that interpretation here.

\subsection{First description}
Let $R(\SU(2))$ be the representation ring of $\SU(2)$, viewed as the ring
of virtual characters. For $m=0,1,2,\ldots$ let $\chi_m\in R(\SU(2))$
be the character of the $m+1$-dimensional irreducible
representation of $\SU(2)$. These form a basis of $R(\SU(2))$ as a
$\Z$-module, and the ring structure is given by
\[\chi_m\chi_{m'}=\chi_{m+m'}+\chi_{m+m'-2}+\cdots+\chi_{|m-m'|}.\]
For $k=0,1,2,\ldots$, the \emph{level $k$ fusion ring} (or
\emph{Verlinde algebra}) is a quotient 
\[ R_k(\SU(2))=R(\SU(2))/I_k(\SU(2))\] 
by the ideal $I_k(\SU(2))$ generated by the character $\chi_{k+1}$.
Additively, the ideal is spanned by the characters
$\chi_{k+1},\,\chi_{2k+3},\,\chi_{3k+5},\ldots$, together with all
characters of the form $\chi_{l'}-(-1)^r\chi_l$ where $l\in
\{0,\ldots,k\}$, and $l'$ is obtained from $l$ by $r$ reflections
across the set of elements $k+1,2k+3,3k+5,\ldots$. It follows that as
an Abelian group, $R_k(\SU(2))$ is free with generators
$\tau_0,\ldots,\tau_k$ the images of $\chi_0,\ldots,\chi_k$. For
example, if $k=4,m=3,m'=4$ we have
\[ 
\chi_3\chi_4=\chi_1+\chi_3+\chi_5+\chi_7\  
  \Rightarrow\ \tau_3\tau_4=\tau_1+\tau_3+0-\tau_3=\tau_1.\]
For any given level $k$, the element $\tau_k\in R_k(\SU(2))$
defines an involution of the group $R_k(\SU(2))$,  
\[ \tau_l\mapsto \tau_l\tau_k=\tau_{k-l}.\]

\subsection{Second description}\label{subsec:two}
Let $q$ be the $2k+4$-th root of unity,
\[ q=e^{\f{i\pi }{k+2}}.\]
Then $I_k(\SU(2))\subset R(\SU(2))$ may be described as the ideal of all characters 
vanishing at all points $j(q^s)$, for  $s=1,\ldots,k+1$. Put
differently, letting
\[ T_{k+2}=\{t\in T|\ t^{2k+4}=e\}\]
be the cyclic subgroup generated by $j(q)$, $I_k(\SU(2))$ is the
vanishing ideal of $T_{k+2}\cap \SU(2)_{\on{reg}}=T_{k+2}^\reg$. 
Hence, for any $t\in T_{k+2}^\reg$ the evaluation map
$\on{ev}_t\colon R(\SU(2))\to \C$ descends to an evaluation map
\[ \on{ev}_t\colon R_k(\SU(2))\to \C,\ \tau\mapsto \tau(t)=\on{ev}_t(\tau).\]
For the basis elements one obtains, by the Weyl character formula,
\[ \tau_l(j(q^s))=\f{q^{(l+1)s}-q^{-(l+1)s}}{q^s-q^{-s}}.\] 
The orthogonality relations 
\begin{equation}\label{eq:orth}
 \sum_{s=1}^{k+1} \f{|q^s-q^{-s}|^2}{2k+4}\,\tau_l(j(q^s)) \tau_{l'}(j(q^s))=\delta_{l,l'}\end{equation}
allow us to recover $\tau\in R_k(\SU(2))$ from the values
$\tau(j(q^s))$ for $s=1,\ldots,k$. The coefficients in this sum may
alternatively be written as
 \[\f{|q^s-q^{-s}|^2}{2k+4}=(\f{k}{2}+1)^{-1}\sin^2(\f{\pi
  s}{k+2}).\] 

\subsection{Third description}
The third way of describing the fusion ring is to write down the
structure constants relative to the basis $\tau_0,\ldots,\tau_k$.  The
\emph{level $k$ fusion coefficient} $N_{l_1,l_2,l_3}^{(k)}$ for $0\le
l_i\le k$ is the multiplicity of $\tau_0$ in the triple product
$\tau_{l_1}\tau_{l_2}\tau_{l_3}$. The fusion coefficients are
invariant under permutations of the $l_i$, and have the
additional symmetry property
$N_{l_1,l_2,l_3}^{(k)}=N_{l_1,k-l_2,k-l_3}^{(k)}$ (coming from
$\tau_{k-l}=\tau_k \tau_{l}$). One has,
\[ \tau_{l_1}\tau_{l_2}=\sum_{l_3=0}^k N_{l_1,l_2,l_3}^{(k)}\tau_{l_3}.\] 
Let $\Delta\subset [0,1]^3$ be the \emph{Jeffrey-Weitsman polytope},
cut out by the inequalities
\[ u_3\le u_1+u_2,\ \ u_1\le u_2+u_3,\ \ u_2\le u_3+u_1,\ \ u_1+u_2+u_3\le
2.\]
Suppose $\Co_i,\ i=1,2,3$ are conjugacy classes of elements
$\exp(u_i\rho)$.  As shown by Jeffrey-Weitsman \cite[Proposition 3.1]{je:bo}, the set
$\{g_1g_2g_3|\ g_i\in\Co_i\}$ contains $e$ if and only if
$(u_1,u_2,u_3)\in \Delta$. Similarly, 
\[ N_{l_1,l_2,l_3}^{(k)}=\begin{cases} 1 & \mbox{if}\ \ \
  l_1+l_2+l_3\ 
  \mbox{even},\ \ (\f{l_1}{k},\f{l_2}{k},\f{l_3}{k})\in
\Delta\\
0 & \mbox{otherwise}
\end{cases}\]

\section{The twisted equivariant $K$-homology of $\SU(2)$}
We will follow the approach to twisted $K$-homology via Dixmier-Douady
bundles.
\subsection{$G$-Dixmier-Douady bundles}
Suppose $G$ is a compact Lie group, acting on a (reasonable)
topological space $X$. A $G$-Dixmier-Douady bundle over $X$ is a
$G$-equivariant bundle $\A\to X$ of $*$-algebras, with typical fiber
$\K(\H)$ the compact operators on a separable Hilbert space $\H$, and
structure group $\on{Aut}(\K(\H))=\on{PU}(\H)$ the projective unitary
group.  Here $\H$ is allowed to be finite-dimensional. A \emph{Morita
isomorphism} between two such bundles $\A_1,\A_2\to X$ is a
$G$-equivariant bundle of $\A_2-\A_1$-bimodules $\E\to X$, such that
$\E$ is locally modeled on the $\K(\H_2)-\K(\H_1)$-bimodule
$\K(\H_1,\H_2)$ of compact operators from $\H_1$ to $\H_2$. We write
\[ \A_1\simeq_\E \A_2.\]
One then also has $\A_2\simeq_{\E^{\on{op}}}\A_1$, where the opposite
  bimodule $\E^{\on{op}}$ is modeled on $\K(\H_2,\H_1)$.  Any two
  Morita isomorphisms $\E,\E'$ between $\A_1,\A_2$ differ by a
  $G$-equivariant line bundle $J$, given as the bundle of bimodule
  homomorphisms:
\[ J=\on{Hom}_{\A_2-\A_1}(\E,\E'),\ \ \ \ \E'=\E\otimes J.\]
Two equivariant Morita isomorphisms $\E,\E'$ will be called
\emph{equivalent} if this line bundle is equivariantly trivial.  By
the Dixmier-Douady theorem \cite{di:ch} (extended to the equivariant
case by Atiyah-Segal \cite{at:twi}), the Morita isomorphism classes of
$G$-Dixmier-Douady bundles $\A\to X$ are classified by an equivariant
\emph{Dixmier-Douady class} $\on{DD}_G(\A)\in H^3_G(X,\Z)$. Put
differently, the Dixmier-Douady class is the obstruction to an
equivariant Morita trivialization $\C\simeq_\E \A$, i.e. an
equivariant Hilbert space bundle $\E$ with an isomorphism
$\A\cong\K(\E)$.
\begin{remark}
  For $G=\{e\}$ the Dixmier-Douady class is realized as a \v{C}ech
  cohomology class, as follows: Choose a cover $\{U_a\}$ of $M$ with
  Morita trivialization $\C\simeq_{\E_a} \A|_{U_a}$. On overlaps, the
  $\E_a$ are related by `transition line bundles',
\[ J_{ab}=\on{Hom}_\A(\E_a,\E_b),\ \ \E_b=\E_a\otimes J_{ab}.\] 
On triple overlaps, one has a trivializing section $\theta_{abc}$ of
$J_{ab}\otimes J_{bc}\otimes J_{ca}$. Taking $U_a$ sufficiently fine,
the $J_{ab}$ are all trivial, and a choice of trivialization makes
$\theta_{abc}$ into a collection of $\U(1)$-valued functions defining
a \v{C}ech cocycle. A different choice of trivialization of the $J_{ab}$
changes the cocycle by a coboundary. The class $\on{DD}(\A)$ equals the
cohomology class of $\theta$, under the isomorphism
$H^2(X,\underline{U(1)})=H^3(X,\Z)$.
\end{remark}

\subsection{The Dixmier-Douady bundle over $\SU(2)$}
We will now give a fairly explicit construction of an equivariant
Dixmier-Douady bundle representing the generator of
$H^3_{\SU(2)}(\SU(2),\Z)=\Z$, using the cover \eqref{eq:cover}.  Let
$\H$ be any $\SU(2)$-Hilbert space, with the property that $\H$ contains all
$T$-weights with infinite multiplicity. (A possible choice is 
$\H=L^2(\SU(2))$ with the left regular representation.) As a
consequence, there exists a $T$-equivariant unitary isomorphism, 
\begin{equation}\label{eq:tisom}
\H\to \H\otimes\C_\rho
\end{equation}
(given by a collection of isomorphisms of the $\mu$-weight spaces with
the $\mu-\rho$-weight spaces). Let
\[ \E_\pm=\SU(2)_\pm\times \H\]
with the diagonal $\SU(2)$-action. By \eqref{eq:useful}, any $\SU(2)$-equivariant
bundle over $\SU(2)_\reg$ is uniquely determined by its restriction to a $T$-equivariant
bundle over $T_{(0,\rho)}$. Let $J\to \SU(2)_\reg$ be the equivariant
line bundle such that $J|_{T_{(0,\rho)}}=T_{(0,\rho)}\times\C_\rho$.
The isomorphism \eqref{eq:tisom} defines a $T$-equivariant isomorphism
\[ \E_-|_{T_{(0,\rho)}}\to  \E_+|_{T_{(0,\rho)}}\otimes J|_{T_{(0,\rho)}}\]
which extends to an $\SU(2)$-equivariant isomorphism 
$\E_-|_{\SU(2)_\reg}\to \E_+|_{\SU(2)_\reg}\otimes J$. This then
defines an isomorphism $\K(\E_-)|_{\SU(2)_\reg}\to \K(\E_+)|_{\SU(2)_\reg}$,
which we use to glue $\K(\E_\pm)$ to a global bundle $\A$.  The bundle
$\A$ represents the generator of $H^3_{\SU(2)}(\SU(2),\Z)=\Z$.  Since
$H^2_{\SU(2)}(\SU(2),\Z)=0$, any other Dixmier-Douady bundle $\A'$
representing the generator is related to $\A$ by a \emph{unique} (up
to equivalence) Morita isomorphism. Again, this can be made quite
explicit: Let $\E_\pm'$ be Morita trivializations of $\A'$, with
transition line bundle $J'$. Then the Morita $\A-\A'$ bimodule is
obtained by gluing $\K(\E'_+,\E_+)$ with $\K(\E'_-,\E_-)$, where the
isomorphism over $\SU(2)_\reg$ is defined by the choice of an
equivariant isomorphism $J'\cong J$ (the latter is unique up to
homotopy).

\subsection{The equivariant Cartan 3-form on $\SU(2)$}\label{subsec:eqCartan}
The equivariant Dixmier-Douady bundle $\A\to \SU(2)$ may be viewed as
a `pre-quantization' of the generator of equivariant Cartan 3-form on
$\SU(2)$. To explain this viewpoint, we need some notation. For any
manifold $M$ with an action of a Lie group $G$, we denote by
$\xi_M\in\mf{X}(M),\ \xi\in\g$ the generating vector fields for the
infinitesimal $\g$-action. That is, $\xi_M(f)=\f{\p}{\p
  u}|_{u=0}(\exp(-u\xi))^*f$ for $f\in C^\infty(M)$.  We let
$(\Om_G^\bullet(M),\d_G)$ denote the complex of equivariant
differential forms
\[ \Om_G^k(M)=\bigoplus_{2i+j=k}(S^i\g^*\otimes\Om^j(M))^G,\]
with equivariant differential
$(\d_G\gamma)(\xi)=\d\gamma(\xi)-\iota(\xi_M)\gamma(\xi)$. For $G$
compact, its cohomology is identified with Borel's equivariant
cohomology $H_G^k(M,\R)$.

Let $\theta^L,\theta^R\in \Om^1(\SU(2),\su(2))$ be the Maurer-Cartan
forms on $\SU(2)$. The \emph{Cartan 3-form} $\eta\in\Om^3(\SU(2))$ is
given in terms of the basic inner product $\cdot$ on $\su(2)$ by
\[ \eta=\f{1}{12}\theta^L\cdot[\theta^L,\theta^L].\]
It is $\d$-closed, and has an equivariantly closed
extension $\eta_{\SU(2)}\in\Om^3_{\SU(2)}(\SU(2))$,
\[ \eta_{\SU(2)}(\xi)=\eta-\hh (\theta^L+\theta^R)\cdot\xi.\]
Let $\varpi\in \Om^2(\su(2))$ be the invariant primitive of
$\exp^*\eta$ defined by the de Rham homotopy operator for the radial
homotopy. The image of the (non-closed) 2-form
$\d\mu-\hh\exp^*(\theta^L+\theta^R)$ under the homotopy operator is 
zero, since its pull-back to any line through the origin vanishes. 
Hence 
\begin{equation}\label{eq:eqdif} 
\exp^*\eta_{\SU(2)}=\d_{\SU(2)}(\varpi-\mu)
\end{equation}
where the
`identity function' $\mu\colon\g\to\g$ is viewed as an element of
$\su(2)^*\otimes \Om^0(\su(2))$. 
\begin{lemma}
For any $G$-manifold with a closed equivariant 3-form
$\gamma\in\Om^3_G(M)$, all $G$-orbits $S\subset M$ acquire
\emph{unique} invariant 2-forms $\om_S\in\Om^2(S)^G$ such that
$\d_G\om_S=i_S^*\gamma$. 
\end{lemma}
The straightforward proof is left to the reader. As special cases,
we obtain 2-forms $\om_\Co$
on the conjugacy classes $\Co\subset \SU(2)$ and $\om_\O$ on the adjoint orbits
$\O\subset \su(2)$ such that 
\[ \d_{\SU(2)}\om_\Co=-\iota_\Co^*\eta_{\SU(2)},\ \ \ 
\d_{\SU(2)}\om_\O=\iota_\O^*(\d\mu).\]
Under the identification of $\su(2)$ with its dual, $\om_\O$ is just the
usual symplectic form on co-adjoint orbits.  Suppose $\Co=\exp(\O)$.
Then \eqref{eq:eqdif} and the uniqueness part of the Lemma imply
\begin{equation}\label{eq:pullb} i_\O^*\varpi=\om_\O-(\exp|_\O)^*\om_\Co.\end{equation}
Let $V\subset \mf{su}(2)$ be the open ball of radius
$\f{1}{\sqrt{2}}$. We have diffeomorphisms
\[ \exp_\pm \colon V\cong \SU(2)_\pm\]
where $\exp_+$ is the restriction of the exponential map, and
$\exp_-=l_c\circ \exp_+$ is its left translate by the central element
$c$.  The inverse maps will be denoted
\[\log_\pm\colon \SU(2)_\pm\to V\subset \su(2).\] 
Let $\varpi_\pm=\log_\pm^*\varpi\in \SU(2)_\pm$. Then
$\d\varpi_\pm=\eta$ over $\SU(2)_\pm$. Furthermore, by Equation
\eqref{eq:eqdif} we have, over $\SU(2)_\pm$,
\begin{equation}\label{eq:eqdiff1}
\d_{\SU(2)}(\varpi_\pm-\log_\pm)=\eta_{\SU(2)}.
\end{equation}
Over $\SU(2)_\reg$, both $\varpi_\pm$ are primitives of $\eta$, hence
their difference is closed. To determine this closed 2-form, 
recall (cf. Equation \eqref{eq:useful}) that $\SU(2)_\reg\cong
T_{(0,\rho)}\times \SU(2)/T$. Let
\[ \Psi\colon \SU(2)_\reg\to \SU(2)/T\]
be the projection to the second factor, and identify $\SU(2)/T$ with the
(co)-adjoint orbit $\O=\SU(2).\rho$.
\begin{lemma}
One has $\varpi_--\varpi_+=\Psi^*\om_\O$ over $\SU(2)_\reg$, where $\O$ is the adjoint
orbit of the element $\rho$.
\end{lemma}
\begin{proof}
By \eqref{eq:eqdiff1} we have 
\[ \d_{\SU(2)}(\varpi_--\varpi_+-(\log_--\log_+))=0\]
over $\SU(2)_\reg$. Thus, $\log_+-\log_-$ serves as a moment map
for the closed invariant 2-form $\varpi_--\varpi_+$.  We claim
\[ \log_+-\log_-=\iota_\O\circ \Psi.\]
Since both sides are ${\SU(2)}$-equivariant, it suffices to compare
the restrictions to $T_{(0,\rho)}\subset \SU(2)_\reg$. Indeed,
$\log_+(\exp(u\rho))=u\rho$ and
$\log_-\exp(u\rho)=\log(\exp(u-1)\rho)=(u-1)\rho$, so the difference
is $(\log_+-\log_-)(\exp(u\rho))=\rho$ as needed. This gives
\[0=\d_{\SU(2)}(\varpi_--\varpi_+ +\iota_\O\circ \Psi)=
\d_{\SU(2)}(\varpi_--\varpi_+ -\Psi^*\om_\O)\]
In particular, $\varpi_--\varpi_+ -\Psi^*\om_\O$ is annihilated by all
contractions with generating vector fields for the conjugation
action.  It is hence enough to show that its pull-back to 
$T_{(0,\rho)}$ is zero. Indeed, by applying the
homotopy operator to $\exp_T^*\iota_T^*\eta_{\SU(2)}=0$, we see that
$\iota_\t^*\varpi=0$, which implies that $\varpi_\pm$ pull back to $0$
on $T$.
\end{proof}
The 2-form $\om_\O$ is the curvature form $\on{curv}(\nabla)$ of the line bundle
$\SU(2)\times_T \C_\rho$, for the unique invariant connection $\nabla$ on
this bundle. Let $J=\Psi^*(\SU(2)\times_T \C_\rho)$ carry the pull-back
connection $\nabla_J$.  The identities
\[ \varpi_--\varpi_+=\on{curv}(\nabla_J),\ \ \ \d\varpi_\pm=\eta\]
say that $(\nabla_J,\varpi_\pm)$ is a `gerbe connection' in the sense of
Chatterjee-Hitchin \cite{ch:co,hi:le}, with $\eta$ as its 3-curvature. Similarly, 
$(\nabla_J,\varpi_\pm-\log_\pm)$ is an equivariant gerbe connection,
with equivariant 3-curvature $\eta_{\SU(2)}$. 

We conclude this Section with an easy proof of the fact that $\eta$
integrates to $1$. Observe that $\partial V=\ol{V}\backslash V$ is the
(co-)adjoint orbit $\O$ of the element $\rho$. It has symplectic
volume $\int_\O\om_\O=1$ by the well-known formula for volume of coadjoint
orbits \cite[Corollary 7.27]{be:he}. Since $\Co:=\exp\O=\{c\}$, we
have $\om_\Co=0$. Hence Equation \eqref{eq:pullb} together with
Stokes' theorem give
\[ \int_{\SU(2)}\eta=\int_{V}\d\varpi=
\int_{\O}\iota_\O^*\varpi=\int_\O\om_\O= 1.\]

\subsection{Twisted $K$-homology}
Let $G$ be a compact Lie group acting on a compact $G$-space $X$.
Given a $G$-Dixmier-Douady bundle $\A\to X$, one defines (following J.
Rosenberg \cite{ros:co}) the twisted $K$-homology group
\[ K_0^G(X,\A)=K^0_G(\Gamma(X,\A)),\]
where the right hand side denotes the $K$-homology group of the
$G-C^*$-algebra of sections of $\A$. (For $K$-homology of
$C^*$-algebras, see \cite{hig:ana,ka:con}.)  The twisted $K$-homology
is a covariant functor: If $\Phi\colon X_1\to X_2$ is an equivariant
map of compact $G$-spaces, together with an equivariant Morita
isomorphism $\A_1\simeq_\E \Phi^*\A_2$, one obtains a push-forward map
\[ \Phi_*\colon K_0^G(X_1,\A_1)\to K_0^G(X_2,\A_2).\]
It is possible to work out
many examples of twisted equivariant $K$-homology groups simply from
its formal properties such as excision, Poincar\'{e} duality and so
on. For $\A=\C$ one obtains
the untwisted $K$-homology groups.  One has a ring isomorphism
\[ K_0^G(\pt)=R(G),\]
where the ring structure on the left hand side is realized
as push-forward under $\pt\times\pt\to \pt$. 
The following is the simplest non-trivial case of the
Freed-Hopkins-Teleman theorem \cite{fr:lo1}. This special case may be proved by an
elementary Mayer-Vietoris argument, see Freed \cite{fr:tw}.
\begin{theorem}
Let $\SU(2)$ act on itself by conjugation, and let $\A\to \SU(2)$ 
be the basic Dixmier-Douady bundle. For all levels $k=0,1,2,\ldots$, 
the $R(\SU(2))$-module homomorphism 
\[ R(SU(2))\cong K_0^{\SU(2)}(\pt)\to K_0^{\SU(2)}(\SU(2),\A^{k+2})\]
given as push-forward under the inclusion of the group unit $\pt\to
\SU(2)$ is onto, with kernel the level $k$ fusion ideal $I_k(\SU(2))$.
It hence defines a ring isomorphism,
\[ R_k(\SU(2))\cong K_0^{\SU(2)}(\SU(2),\A^{k+2}).\]
\end{theorem}

\subsection{The $K$-homology fundamental class}\label{subsec:funclass}
Recall that for $n$ even, the complex Clifford algebra
$\Cl(n)=\Cl(\R^n)$ admits a unique (up to isomorphism) irreducible
$*$-representation. Concretely, the identification $\R^n\cong
\C^{n/2}$ gives a Clifford action on the standard spinor module
$\mathsf{S}=\wedge \C^{n/2}$. This realizes the Clifford algebra as a 
matrix algebra, $\Cl(n)=\on{End}(\mathsf{S})$. Given $A\in \SO(n)$
there exists a unitary transformation $U\in \on{U}(\mathsf{S})$,
unique up to a scalar, such that $A(v).U(z)=U(v.z)$ for $v\in \R^n,\ 
z\in \mathsf{S}$.  The set of such \emph{implementers} $U$ forms a
closed subgroup of $\on{U}(\mathsf{S})$, denoted $\Spin_c(n)$, and the map
taking $U$ to $A$ makes this group into a central extension
\[ 1\to \U(1)\to \Spin_c(n)\to \SO(n)\to 1.\]
If $M$ is an oriented Riemannian $G$-manifold of even dimension $n$,
then its Clifford algebra bundle $\Cl(TM)$ is a $G$-equivariant bundle
of complex matrix algebras. It is thus a $G$-Dixmier-Douady bundle.
Its Dixmier-Douady class is the third integral equivariant
\footnote{We remark that for $G$ compact and simply connected, the
  vanishing of $W^3_G(M)$ is equivalent to the vanishing of the
  non-equivariant Stiefel-Whitney class $W^3(M)$, since the map
  $H^3_{G}(M,\Z)\to H^3(M,\Z)$ is injective (cf. \cite{kre:pr}).}
Stiefel-Whitney class, $W^3_G(M)\in H^3_G(M,\Z)$.  As pointed out by
Connes \cite{con:non} and Plymen \cite{ply:st}, an equivariant
$\Spin_c$-structure on $M$ is exactly the same thing as an equivariant
Morita trivialization of $\Cl(TM)$. Indeed, given an equivariant lift
$P_{\Spin_c}(M)\to P_{SO}(M)$ of the $\SO(n)$-frame bundle to
the group $\Spin_c(n)$, the Morita trivialization is defined by the
bundle of spinors $\S=P_{\Spin_c}(M)\times_{\Spin_c(n)}\mathsf{S}$.
Conversely, given an equivariant Morita trivialization
$\Cl(TM)\simeq_{\S}\C$, on obtains a lift of the structure group: The
fiber of the bundle $P_{\Spin_c}(M)$ at $m\in M$ is the set of pairs
$(A,U)$, where $A\colon T_mM\to \R^n$ is an oriented orthonormal
frame, and $U\colon \S_m\to \mathsf{S}$ is a unitary isomorphism
intertwining the Clifford actions of $v\in T_mM$ and $A(v)\in \R^n$.

The Clifford bundle $\Cl(TM)$ is naturally a $\Cl(TM)-\Cl(TM)$ bimodule. 
Using the canonical anti-automorphism of $\Cl(TM)$, it may also be viewed
as a module over $\Cl(TM)\otimes\Cl(TM)$, defining a Morita 
trivialization of the latter. Given any $\Spin_c$-structure $\S$, 
one obtains a Hermitian line bundle 
\[ \L:=\L(\S)=\on{Hom}_{\Cl(TM)\otimes \Cl(TM)}(\Cl(TM),\S\otimes\S)\]
called the \emph{$\Spin_c$-line bundle}.  Twisting $\S$ by a line
bundle $L$ changes the $\Spin_c$-line bundle as follows,
\[ \L(\S\otimes L)=\L(\S)\otimes L^2.\] 
For any equivariant $\Spin_c$-structure on an even-dimensional
manifold, the class of the $\Spin_c$-Dirac operator defines a
\emph{fundamental class} in equivariant $K$-homology. In the absence
of a $\Spin_c$-structure, there is still a fundamental class, but as
an element
\[ [M]\in K_0^G(M,\Cl(TM))\] 
in twisted
$K$-homology. \footnote{More precisely, one has to view $\Cl(TM)$ as a
  $\Z_2$-graded Dixmier-Douady bundle, and work with the twisted
  $K$-homology for such $\Z_2$-graded bundles.}
For an explicit construction of $[M]$, see Kasparov \cite{ka:con}.
Below, we will construct elements of
$R_k(\SU(2))=K_0^{\SU(2)}(\SU(2),\A^{k+2})$ as push-forwards of $[M]$
under $\SU(2)$-equivariant maps $\Phi\colon M\to \SU(2)$.  In order to
define such a push-forward, we need an equivariant Morita isomorphism
\[\Cl(TM)\simeq_\E \Phi^*\A^{k+2}.\]
We will explain how such a `twisted $\Spin_c$-structure' arises for
pre-quantized q-Hamiltonian $\SU(2)$-spaces. The counterpart to the
$\Spin_c$-line bundle is the Morita isomorphism
$\Phi^*\A^{2k+4}\simeq_{\mathcal{K}} \C$ given by
\[ \mathcal{K}=\on{Hom}_{\Cl(TM)\otimes \Cl(TM)}(\Cl(TM),(\E\otimes\E)^{\on{op}}).\]

\section{q-Hamiltonian $\SU(2)$-spaces}
\subsection{Basic definitions}
Let $G$ be a compact Lie group, with Lie algebra $\g$. Given an
invariant inner product $B$ on its Lie algebra, define the equivariant
Cartan 3-form
\[ \eta^{(B)}_G(\xi)=\f{1}{12}B(\theta^L,[\theta^L,\theta^L])-\hh B(\theta^L+\theta^R,\xi).\]
A \emph{q-Hamiltonian $G$-space} (relative to the inner product $B$) is a triple
$(M,\om,\Phi)$ where $M$ is a $G$-manifold, $\omega$ is an invariant
2-form, and $\Phi\colon M\to G$ an equivariant smooth map, called the 
\emph{moment map}, such that 
\begin{itemize}
\item[(i)] $\d_G\om=-\Phi^*\eta_G^{(B)}$, 
\item[(ii)] $\ker\omega\cap \ker(\d\Phi)=0$ everywhere.
\end{itemize}
\begin{remark}
  If $G=T$ is a torus, this is just the usual definition of a
  symplectic $T$-space with torus-valued moment map. Indeed, Condition
  (i) in this case says $\d\om=0$ and
  $\om_m(\xi_M(m),v)=-B(\theta_T(\d_m\Phi(v)),\xi)$ for all
  $\xi\in\g,\ v\in T_mM$. Hence it implies $\ker(\om)\subset
  \ker(\d\Phi)$, whence (ii) simplifies to $\ker(\om)=\{0\}$.
  For general $G$, a similar argument shows that $\ker(\omega_m)$ is
  spanned by all $\xi_M(m)$ such that $\Ad_{\Phi(m)}\xi+\xi=0$.
\end{remark}
Basic examples of q-Hamiltonian $G$-spaces are the conjugacy classes
$\Co\subset G$, with moment map the embedding. The \emph{double}
$D(G)=G\times G$, with $G$ acting by conjugation and with moment
$\Phi(a,b)=aba^{-1}b^{-1}$, is another example. The 2-form is,
\[ \om=\hh a^*\theta^L\cdot b^*\theta^R+\hh a^*\theta^R\cdot
b^*\theta^L+
\hh (ab)^*\theta^L\cdot (a^{-1}b^{-1})^*\theta^R,\]
where, for example, $a^{-1}b^{-1}$ denotes the map $(a,b)\mapsto
a^{-1}b^{-1}$.  If $G'$ is the quotient of $G$ by a finite subgroup of
$Z(G)$, then the moment map, action and 2-form on $D(G)$ descends to
$D(G')$, so that $D(G')$ is again a q-Hamiltonian $G$-space.

Given two q-Hamiltonian $G$-spaces
$(M_i,\om_i,\Phi_i),\ i=1,2$, their product $M_1\times M_2$ with the
diagonal $G$-action, moment map $\Phi_1\Phi_2$, and 2-form
$\om_1+\om_2+\hh B(\Phi_1^*\theta^L,\Phi_2^*\theta^R)$ is again a
q-Hamiltonian $G$-space. This is called the \emph{fusion product} of
$M_1,M_2$.  The \emph{symplectic quotient} of a q-Hamiltonian
$G$-space is $M\qu G=\Phinv(e)/G$. Similar to the Hamiltonian theory,
$e$ is a regular value of $\Phi$ if and only if $G$ acts locally
freely on $\Phinv(e)$, and in this case $M\qu G$ is a symplectic
orbifold. (If $e$ is a singular value, then $M\qu G$ is a singular
symplectic space as defined in \cite{sj:st}.) More generally, given a
conjugacy class $\Co$ one can define a symplectic quotient
\[ M\qu_\Co G = (M\times\Co)\qu G.\]
It was shown in \cite{al:mom} that moduli spaces of flat $G$-bundles
over compact oriented surfaces $\Sigma_h^r$ of genus $h$ with $r$
boundary circles, with boundary holonomies in prescribed conjugacy
classes $\Co_j$, are symplectic quotients
\[ M(\Sig_h^r,\Co_1,\ldots,\Co_r)=(\underbrace{D(G)\times\cdots\times
D(G)}_{{h \ \on{times}}}\times\Co_1\times\cdots\Co_r)\qu G.\]
We now specialize to $q$-Hamiltonian $\SU(2)$-spaces $(M,\om,\Phi)$, with $B$ the
basic inner product. Put $M_\pm=\Phinv(\SU(2)_\pm)$, and let 
\[ \begin{split}
\om_{0,\pm}&=\om+\Phi^*\varpi_\pm,\\ 
\Phi_{0,\pm}&=\log_\pm\circ \Phi.\end{split}\]
Then
\[ \d_{\SU(2)}(\om_{0,\pm}-\Phi_{0,\pm})=
\d_{\SU(2)}(\om+\Phi^*(\varpi_\pm-\log_\pm))=0.\]
That is, $\om_{0,\pm}$ is closed, with $\Phi_{0,\pm}$ as a moment map.
Using condition (ii) above one can show \cite{al:mom} that
$\om_{0,\pm}$ are non-degenerate, i.e.  \emph{symplectic}.  Thus,
$(M_\pm,\om_{0,\pm},\Phi_{0,\pm})$ are ordinary (symplectic)
Hamiltonian $\SU(2)$-spaces. In particular, $M_\pm$ are
even-dimensional, with a natural orientation. If $M$ is compact and
connected, then the spaces $M_\pm$ are connected. (This follows from
the convexity properties and the fiber connectivity of 
group-valued moment maps \cite{al:mom}.)

Conversely, $(M,\om,\Phi)$ is determined by the pair of Hamiltonian
$\SU(2)$-spaces $(M_\pm,\om_{0,\pm},\Phi_{0,\pm})$.  This
correspondence reduces many properties of q-Hamiltonian spaces to
standard facts about ordinary Hamiltonian spaces. It is also used to
construct q-Hamiltonian spaces, as in the following example.
\subsection{Example: The 4-sphere}
The following construction of a q-Hamiltonian structure of $S^4$ is
taken from \cite{al:du}.  An independent construction due to
Hurtubise-Jeffrey \cite{hu:re1} was later generalized by
Hurtubise-Jeffrey-Sjamaar \cite{hu:imp} to define the structure of a
q-Hamiltonian $\SU(n)$-space on $S^{2n}$, for any $n$.

Let $\C^2$ carry the standard $\SU(2)$-action and the standard 
symplectic structure $\om_0=\f{i}{2}(\d z_1 \wedge \d\ol{z}_1+\d
z_2\wedge\d\ol{z}_2)$. The moment map for the $\SU(2)$-action can be
written, for $z\not=0$, as
\[\Phi_0(z)=-i\pi^2||z||^2 P(z)+i\pi^2 ||z||^2 (I-P(z)),\]
where $P(z)$ is the projection operator, \
\[ P(z)=||z||^{-2}\,
\left(\begin{array}{c}{z_1}\\{z_2}\end{array}\right)\left(\begin{array}{c}z_1\\z_2\end{array}\right)^\dagger=\f{1}{||z||^2}\left(\begin{array}{cc}
    |z_1|^2 & z_1\ol{z}_2\\
    \ol{z}_1 z_2& |z_2|^2
\end{array}\right).\]
Hence, 
\[ \exp(\Phi_0(z))=e^{-i\pi^2||z||^2} P(z)+e^{i\pi^2||z||^2}(I-P(z)).\]
Let $V\subset \su(2)$ be the open ball of radius $\f{1}{\sqrt{2}}$
(cf. Section \ref{subsec:eqCartan}). 
We have $||\Phi_0(z)||=\f{1}{\sqrt{2}}\pi ||z||^2$, so that
\[ S^4_\pm:= \Phi_0^{-1}(V)=\{z\in \C^2|\ \pi ||z||^2<1\}.\] 
Define a diffeomorphism $F$ of the annulus $0< \pi ||z||^2<1$ by 
\[F(z_1,z_2)=(-\ol{z_2},\ol{z}_1)\ \sqrt{\textstyle{\f{1}{\pi ||z||^2}}-1}.\]
Then $F$ is equivariant, with $\pi||F(z)||^2=1-\pi||z||^2$.  Gluing
the charts $S^4_\pm$ under $F$ one obtains a 4-sphere $S^4$ with an
action of $\SU(2)$.

Put $\Phi_+=\exp\Phi_0$ and  $\Phi_-=l_c\circ \exp\Phi_0=-\exp\Phi_0$.  
The diffeomorphism $F$ satisfies $P(F(z))=I-P(z)$, and therefore,
\[ \Phi_+(F(z))=\exp(\Phi_0(F(z)))=-\exp(\Phi_0(z))=\Phi_-(z).\]
Hence $\Phi_\pm$ glue to a global equivariant map $\Phi\colon S^4\to
\SU(2)$. Similarly, the 2-forms $\om_\pm=\om_{0}+\Phi_0^*\varpi$ glue
\footnote{To check that these 2-forms agree on the overlap
  $S^4_\reg=S^4_+\cap S^4_-$, it suffices to consider their pull-back
  to symplectic cross-sections as in Section \ref{sec:cross}.}
to a global invariant 2-form $\om\in \Om^2(S^4)$, defining a
q-Hamiltonian $\SU(2)$-space $(S^4,\om,\Phi)$.
\begin{remark}
  The space $S^4$ carries an involution $I\colon S^4\to S^4$, given in
  charts by the complex conjugation. It has the equivariance property
  $I(g.x)=I(g).I(x)$ relative to the involution of $\SU(2)$ given by
  complex conjugation of matrices, $I(A)=\ol{A}$. The involution
  satisfies, $I^*\om=-\om$ and $I^*\Phi=\ol{\Phi}$. The fixed point
  set of the involution is a 2-sphere $S^2\subset S^4$. The theory of
  anti-involutions of $q$-Hamiltonian $G$-spaces was developed in
  recent work of Schaffhauser \cite{sch:re}, who established an
  analogue of the convexity results of Duistermaat \cite{du:co} and
  O'Shea-Sjamaar \cite{os:th} in this context.
\end{remark}

\begin{remark}
  It is well-known that the complement of the zero section in
  $T^*(S^2)$ is $\SU(2)$-equivariantly symplectomorphic to the
  complement of the origin in $\C^2$. One may thus modify the
  construction above, and obtain examples where the the fiber over $e$
  or over $c$ (or both) is a 2-sphere rather than a point. The four
  examples obtained in this way are the complete list of $4$-dimensional
  q-Hamiltonian $\SU(2)$-spaces with surjective moment map.
\end{remark}

\section{Cross-sections}\label{sec:cross}
Let $(M,\om,\Phi)$ be a q-Hamiltonian $\SU(2)$-space.  By the
q-Hamiltonian cross-section theorem \cite{al:mom}, the pre-image
\begin{equation}\label{eq:crosssection} Y=\Phinv(T_{(0,\rho)})\end{equation}
is a q-Hamiltonian $T$-space $(Y,\om_Y,\Phi|_Y)$, with 2-form
$\om_Y=i_Y^*\om$. In particular, $\om_Y$ is symplectic. Letting
$\Phi_Y\colon Y\to (0,\rho)\subset\t$ with $\exp\Phi_Y=\Phi|_Y$, it is
immediate that $(Y,\om_Y,\Phi_Y)$ is an ordinary Hamiltonian
$T$-space. We have,
\[ M_\reg=M_+\cap M_-=\SU(2)\times_T Y\]
and 
\[ TM|_Y=TY\oplus \t^\perp,\] 
where the second summand is embedded by the generating vector fields.
This splitting is $\om$-orthogonal, and the 2-form on
$Y\times\t^\perp$ is given at $y\in Y$, with $g=\Phi(y)\in
T_{(0,\rho)}$, by $(\xi_1,\xi_2)\mapsto \hh
((\Ad_g-\Ad_{g^{-1}})\xi_1,\xi_2)$. Note that since the pull-back of
$\varpi_\pm$ to $T_{(0,\rho)}$ is zero, the 2-forms $\om_{0,\pm}$ both
pull back to $\om_Y$. Similarly 
\[ \Phi_{0,+}|_Y =\Phi_Y=\Phi_{0,-}|_Y+\rho.\]
That is, $(Y,\om_Y,\Phi_Y)$ may also be viewed as symplectic
cross-section of $M_\pm$. (To be precise, in the case of $M_-$, it is
the \emph{opposite} cross-section, given as the pre-image of
$(-\infty,0)\subset\t$ under $\Phi_{0,-}$.) The 2-forms on the bundles 
$Y\times \t^\perp$ induced by $\om_{0,\pm}$ are, 
\[(\xi_1,\xi_2)\mapsto \ad_{\mu_\pm}\xi_1\,\cdot\xi_2,\] 
where $\mu_+=\Phi_{0,+}(y)$ and $\mu_-=\Phi_{0,-}(y)$.

The space $Y$ is only a `partial' cross-section for $M$, since it
leaves out the subsets $\Phinv(e),\ \Phinv(c)$. On the other hand, the
`full' cross-section $\ti{Y}=\Phinv(T_{[0,\rho]})$ is usually not a
manifold, let alone symplectic. However, following
Hurtubise-Jeffrey-Sjamaar \cite{hu:imp} one can `implode' $\ti{Y}$ to
obtain a symplectic $T$-space $X$, which is a symplectic orbifold
under regularity conditions. As a topological space, the
\emph{imploded cross-section} is a quotient space
\[ X=\Phinv(T_{[0,\rho]})/\sim,\]
where the equivalence relation divides out the $\SU(2)$-action on both
$\Phinv(e)$ and on $\Phinv(c)$. We have a
decomposition of $X$ into three symplectic spaces,
\begin{equation}\label{eq:decomp} X=(M\qu \SU(2)) \cup Y \cup (M\qu_c \SU(2))\end{equation}
The action of $T\subset \SU(2)$ on $\Phinv(T_{[0,\rho]})$ descends to
an action on $X$, and the map
$\Phinv(T_{[0,\rho]})\to [0,\rho]\subset \t$ descends to a
$T$-equivariant map
\[ \Phi_X\colon X\to \t.\]
Let
\[ X_+=(M\qu \SU(2)) \cup Y,\ \
X_-=Y \cup (M\qu_c \SU(2)),\]
so that $X_\pm$ are the imploded cross-sections of $M_\pm$.  View
$M_\pm$ as Hamiltonian $\SU(2)$-spaces with 2-forms $\om_{0,\pm}$, and
let $\C^2$ carry the standard structure as a Hamiltonian
$\SU(2)$-space.
\begin{proposition}
Suppose $\SU(2)$ acts locally freely (resp. freely) on
$\Phinv(e),\Phinv(c)$. Then the imploded cross-section $X$ admits a unique
structure of a symplectic orbifold (resp. symplectic manifold), 
such that the open subsets $X_\pm$ are symplectic quotients, 
\[ X_\pm=(M_\pm\times\C^2)\qu \SU(2).\]
Furthermore, 
\begin{enumerate}
\item
The restriction of $\Phi_X$ to $X_\pm$ is smooth, and is a moment map
for the action of $T\cong \U(1)$. 
\item 
The Hamiltonian $T$-space $(Y,\om_Y,\Phi_Y)$ is embedded as an open symplectic submanifold of $X$.
\item 
$M\qu \SU(2)$ is a symplectic suborbifold (resp. submanifold), with normal bundle
$\Phinv(e)\times_{\SU(2)}\C^2$. The $\U(1)$ action on the normal
bundle is with weights $(-1,-1)$.
\item 
$M\qu_c \SU(2)$ is a symplectic suborbifold (resp. submanifold), with normal bundle
$\Phinv(c)\times_{\SU(2)}\C^2$.  The $\U(1)$-action on the normal bundle
is with weights $(1,1)$.
\end{enumerate}
\end{proposition}
Thus, $X$ is obtained by gluing the Hamiltonian imploded
cross-sections for $(M_\pm,\om_{0,\pm},\Phi_{0,\pm})$. For the case
$G=\SU(2)$, the imploded cross-sections construction was introduced by
Eugene Lerman as an $\SU(2)$-counterpart of symplectic cutting. Its
basis properties for Hamiltonian $\SU(2)$-spaces are described in
\cite[Appendix]{me:sym}, and directly imply imply the properties for
q-Hamiltonian $\SU(2)$-spaces.
%
\begin{remark}
  More intrinsically, the imploded cross-section can directly be
  constructed as a q-Hamiltonian symplectic quotient $X=(M\times
  S^4)\qu \SU(2)$. This is the approach taken in
  \cite{hu:re1,hu:imp}. However, in this paper we will have more use
  for the construction in terms of ordinary Hamiltonian quotients.
\end{remark}

\section{The canonical `twisted $\Spin_c$-structure'}
Choose invariant almost complex structures on $M_\pm$, which are
compatible with $\om_{0,\pm}$ in the sense that each tangent space is
isomorphic to $\C^{n/2}$ with the standard complex structure and standard
symplectic form. The almost complex structure defines spinor modules
\[ \S_{0,\pm}=\wedge_\C TM_\pm\to M_\pm\] 
for the Clifford bundles $\Cl(TM)|_{M_\pm}$, where the notation
$\wedge_\C$ denotes the complex exterior powers of $TM_\pm$ relative
to the given complex structure. On the overlap $M_+\cap M_-=M_\reg$, the two spinor bundles
differ by $\on{Hom}_{\Cl(TM)}(\S_{0,+},\S_{0,-})$.

\begin{proposition}
The line bundle $\on{Hom}_{\Cl(TM)}(\S_{0,+},\S_{0,-})$ is equivariantly
isomorphic to the pull-back $\Phi^*(J^{\otimes 2})$. 
\end{proposition}
\begin{proof}
  An $\SU(2)$-invariant almost complex structure on
  $M_\reg=\SU(2)\times_T Y$ is equivalent to a $T$-invariant complex
  structure on the bundle $TM|_Y=TY\oplus \t^\perp$.  This bundle
  carries two symplectic structures, defined by the 2-forms
  $\om_{0,\pm}$ on $M_\pm$. Pick a $T$-invariant compatible structure
  on the bundle $TY$. Its sum with the complex structure on
  $\t^\perp$, coming from the identification $\t^\perp\cong
  \C_\alpha$, is compatible with $\om_{0,+}$. Similarly its sum with
  the complex structure on $\t^\perp$, coming from the identification
  $\t^\perp\cong \C_{-\alpha}$, is compatible with $\om_{0,-}$.  The
  corresponding spinor bundles $\ti{\S}_{0,\pm}|_Y\to Y$ are related by a
  twist by a $T$-equivariant line bundle, corresponding to the change
  of the complex structure on $\t^\perp$ to its opposite. Clearly,
  this is the line bundle $Y\times \C_\alpha=Y\times (\C_\rho)^2$: 
\[ \ti{\S}_{0,-}|_Y=\ti{\S}_{0,+}|_Y\otimes (Y\times (\C_\rho)^2).\]
Extending to $M_{reg}$, and using the definition of $J\to \SU(2)_\reg$ we obtain 
\[ \ti{\S}_{0,-}=\ti{\S}_{0,+}\otimes \Phi^* J^2.\]
But $\ti{S}_{0,\pm}$ are equivariantly isotopic to $\S_{0,\pm}$, since any two
choices of equivariant compatible almost complex structures are
isotopic. Hence we also have $\S_{0,-}\cong \S_{0,+}\otimes \Phi^* J^2$, or
equivalently $\on{Hom}_{\Cl(TM)}(\S_{0,+},\S_{0,-})\cong \Phi^* J^2$.
\end{proof}

Equivalently, we can express this result as follows: 
\begin{proposition}\label{prop:spinc}
For any q-Hamiltonian $\SU(2)$-space $(M,\om,\Phi)$, there is a
distinguished (up to equivalence) $\SU(2)$-equivariant Morita isomorphism 
\begin{equation}\label{eq:cansp}
\Phi^*\A^{2}\simeq_\S \Cl(TM),
\end{equation}
\end{proposition}
\begin{proof}
  Let $\F_\pm\to \SU(2)_\pm$ define Morita trivializations $\C
  \simeq_{\F_\pm} \A^{2}$. Fix isomorphisms $\F_-\cong \F_+\otimes
J^2$ and $\S_{0,-}\cong \S_{0,+}\otimes \Phi^*J^2$ on intersections. 
The desired Morita $\Cl(TM)-\Phi^*\A^2$ bimodule $\S$ is then obtained by gluing
the bundles $\S_\pm=\on{Hom}_\C(\Phi^*\F_\pm,\S_{0,\pm})$, using that 
\[ \on{Hom}_\C(\Phi^*\F_-,\ \S_{0,-})\cong 
\on{Hom}_\C(\Phi^*(\F_+\otimes
J^2),\ \ \S_{0,+}\otimes \Phi^*J^2)=\on{Hom}_\C(\Phi^*\F_+,\ \S_{0,+})\]
on the intersection. 
\end{proof}

We refer to the Morita isomorphism \eqref{eq:cansp} as the canonical
twisted $\Spin_c$-structure of a q-Hamiltonian manifold. 

\begin{remark}
In particular, we see that the third integral Stiefel-Whitney class of
any q-Hamiltonian $\SU(2)$-space satisfies 
\[ W^3(M)=2\Phi^* x\]
where $x\in H^3(\SU(2),\Z)$ is the generator. Since this is a
2-torsion class, it follows that $4\Phi^*x=0$. The fact that $\Phi^*x$
is torsion is a consequence of the condition $\d\om=-\Phi^*\eta$. The
more precise statement relies on the minimal degeneracy condition
$\ker(\om)\cap \ker(\d\Phi)=0$.
\end{remark}

\section{Pre-quantization of q-Hamiltonian $\SU(2)$-spaces}
Suppose $(M,\om,\Phi)$ is a q-Hamiltonian $\SU(2)$-space. The
conditions $\d\om=-\Phi^*\eta$ and $\d\eta=0$ mean that the pair
$(\om,-\eta)$ defines a cocycle for the relative de Rham complex
\footnote{Recall that for any morphism of cochain complexes
$F^\bullet\colon C^\bullet\to \ti{C}^\bullet$, the relative cohomology
$H^\bullet(F)$ is the cohomology of the algebraic mapping cone
$(\ti{C}^{k-1}\oplus C^k,\d)$, with differential $\d(x,y)=(F(y)-\d
x,\d y)$. In our case $F=\Phi^*$, acting on differential forms or on
singular cochains, and we write $H(\Phi,\cdot)$ for the relative
cohomology.}  $\Om^\bullet(\Phi)$. For $k>0$, we define a \emph{level
$k$ pre-quantization} of $(M,\om,\Phi)$ to be a lift of the class
$k[(\om,-\eta)]\in H^3(\Phi,\R)$ to a class in $H^3(\Phi,\Z)$.
\begin{remark}
  One can similarly define an \emph{equivariant} level $k$
  pre-quantization to be an integral lift of
  $k[(\om,-\eta_{\SU(2)})]\in H^3_{\SU(2)}(\Phi,\R)$. However, the
  equivariance is automatic: Indeed, for any simply connected compact
  Lie group $G$, and any $G$-space $M$ one has $H^p_G(M,\Z)=H^p(M,\Z)$
  for $p\le 2$, and if $\Phi\colon M\to G$ is an equivariant map one
  has $H^p_G(\Phi,\Z)=H^p(\Phi,\Z)$ for $p\le 3$. See e.g.
  \cite{kre:pr}.
\end{remark}

\begin{lemma}
  If $(M,\om,\Phi)$ admits a level $k$ pre-quantization, then the set
  of such pre-quantizations is a principal homogeneous space under the
  group $\on{Tor}(H^2(M,\Z))$ of flat line bundles over $M$.
\end{lemma}
\begin{proof}
  Clearly, the set of pre-quantizations is a principal homogeneous
  space under $\on{Tor}(H^3(\Phi,\Z))$. Since $H^3(\SU(2),\Z)=\Z$ has no
  torsion, $\on{Tor}(H^3(\Phi,\Z))$ lies in 
  the image of the map $H^2(M,\Z)\to H^3(\Phi,\Z)$ in the long
  exact sequence for relative cohomology. But this map is injective since
  $H^2(\SU(2),\Z)=0$, and hence restricts to an isomorphism of the
  torsion subgroups.
\end{proof}

The class $k[(\om,-\eta)]$ is integral if and only if it takes integer
values on all relative 3-cycles: That is, for every smooth singular
2-cycle $\Sigma\in C_2(M)$, and every smooth singular $3$-chain $\Gamma\in
C_3(\SU(2))$ bounding $\Phi(\Sigma)$, we must have
\begin{equation}\label{eq:criterion1} k(\int_\Gamma\eta+\int_\Sigma\om)\in\Z.
\end{equation}
(Given $\Sigma$, it is actually enough, by the integrality of $\eta$,
to check the condition for for \emph{some} $\Gamma$ bounding
$\Phi(\Sigma)$.) If $H^2(M,\R)=0$, there is a much simpler criterion
\cite{kre:pr}: Let $x\in H^3(\SU(2),\Z)$ be the generator. Since
$\Phi^*[\eta]=0$, the class $\Phi^*x$ is torsion. If $H^2(M,\R)=0$,
then $(M,\om,\Phi)$ is pre-quantizable at level $k$ if and only if
\begin{equation}\label{eq:criterion2}
k\Phi^*x=0.\end{equation}

\begin{proposition}
The conjugacy class $\Co$ of $t\in T_{[0,\rho]}\subset\SU(2)$ 
is pre-quantizable at level $k$ if and only if $t=\exp(\f{n}{k} \rho)$ for some 
$n\in \{0,1,\ldots,k\}$. 
\end{proposition}
\begin{proof}
  It is enough to check Criterion \eqref{eq:criterion1} for $\Sigma=\Co$. Write
  $t=\exp(u\rho)$ with $u\in [0,1]$. Let
  $\O$ be the adjoint orbit of $u\rho$, so that $\Co=\Phi(\O)$. As above, let
  $V\subset \su(2)$ be the open ball of radius $\f{1}{\sqrt{2}}$.
  Then $\O$ is the boundary of $V_u=uV$, and we compute,
    with $\Gamma=\Phi(\ol{V_u})$,
\[\int_\Gamma\eta=\int_{V_u}\exp^*\eta=\int_{V_u}\d\varpi=\int_\O
i_\O^*\varpi=\int_\O\om_\O-\int_\Co\om_\Co.\]
Hence
\[ k(\int_\Gamma\eta+\int_\Co\om_\Co)=k\int_\O\om_\O\]
which is an integer if and only if the orbit through $ku\rho$ 
is integral, i.e. $ku\in\Z$. 
\end{proof}

\begin{proposition}
  The 4-sphere $S^4$ and the double $D(\SU(2))$ are pre-quantizable at
  any integer level $k$. More generally, this is the case for any
  q-Hamiltonian $\SU(2)$-space $(M,\om,\Phi)$ with vanishing second
  homology. The double $D(\SO(3))$ (viewed as a q-Hamiltonian
  $\SU(2)$-space) is pre-quantizable at level $k$ if and only if $k$
  is even. 
\end{proposition}
The condition for $D(\SO(3))$ was first obtained by Derek Krepski
\cite{kre:pr}.
\begin{proof}
  In each of these examples we have $H^2(M,\R)=0$, hence it suffices
  to find all $k$ such that $k\Phi^*x=0$. For $M=S^4$, one has
  $\Phi^*x=0$ since $H^3(S^4,\Z)=0$. For $M=D(\SU(2))$, one again has
  $\Phi^*x=0$, by the properties of $x$ under group multiplication and
  inversion ($\on{Mult}^*x=\pr_1^* x+\pr_2^*x,\ \on{Inv}^*x=-x$.) For
  $M=D(\SO(3))$, one checks that the torsion subgroup of $H^3(M,\Z)$
  is $\Z_2$, so that $M$ is pre-quantizable at either all levels or at
  all even levels. We claim that $M$ is \emph{not} pre-quantizable at
  level $1$.  To see this consider the symplectic submanifold
  $T'\times T'\subset D(\SO(3))$, where $T'$ is the maximal torus in
  $\SO(3)$ given as the image of $T$. For the symplectic volume one
  finds, (see \ref{subsec:double} below)
\[ \on{vol}(T'\times T')=\f{1}{4} \on{vol}(T\times T)=\f{2}{4}=\f{1}{2}.\]
  By Criterion \eqref{eq:criterion1}, with $\Sig=T'\times T'$ and
  $\Gamma=\emptyset$, the pre-quantized levels $k$ must satisfy
  $k\int_\Sigma\om\in\Z$, hence they must be even.
\end{proof}
Finally, we remark that if $(M_i,\om_i,\Phi_i)$ are pre-quantized at
level $k$, then their fusion product $M_1\times M_2$ inherits a
pre-quantization at level $k$. 

For an ordinary Hamiltonian $\SU(2)$-space $(M,\om_0,\Phi_0)$, a
pre-quantization is an integral lift of the class of the equivariant
symplectic form. More generally, by a \emph{level $k$ pre-quantization}
of such a space we mean a pre-quantization of $(M,k\om_0,k\Phi_0)$.
Geometrically, the lift is realized as the equivariant Chern class of
an equivariant \emph{pre-quantum line bundle} over $M$.
\begin{proposition}
A level $k$ pre-quantization of a  
q-Hamiltonian $\SU(2)$-space $(M,\om,\Phi)$ is
equivalent to a pair of level $k$ pre-quantizations of 
the  Hamiltonian $\SU(2)$-spaces
$(M_\pm,\ \om_{0,\pm},\ \Phi_{0,+})$, with the property that the 
the  pre-quantum line bundles $L_\pm\to M_\pm$ satisfy
  \[ L_-\cong L_+\otimes \Phi^*J^k\] 
on the overlap $M_\reg=M_+\cap M_-$.
\end{proposition}
\begin{proof}
  Let $\Phi_\pm\colon M_\pm \to \SU(2)_\pm$ be the restrictions of
  $\Phi$. Since $\SU(2)_+,\SU(2)_-$ retract onto $e,c$ respectively,
  the long exact sequences in relative cohomology give isomorphisms
  $H^2(M_\pm,\cdot)\xrightarrow{\cong} H^3(\Phi_\pm,\cdot)$, and a
  commutative diagram,
\[ \begin{CD} H^3(\Phi,\Z) @>>> H^3(\Phi_\pm,\Z) \cong H^2(M_\pm,\Z)\\
@VVV @VVV \\ H^3(\Phi,\R) @>>> H^3(\Phi_\pm,\R) \cong H^2(M_\pm,\R)
\end{CD}\]
The lower horizontal map is given on $k[(\om,-\eta)]$ by
\[ k[(\om,-\eta)] \mapsto k[\om_\pm+\Phi_\pm^*\varpi_\pm]=k[\om_{0,\pm}].\]
To give a parallel discussion of the upper horizontal map, let
$C^k(\cdot,R)=\on{Hom}(C_k(\cdot),R)$ denote the complex of smooth
singular cochains, with coefficient in the ring $R$.  We have two
natural cochain maps,
\[ C^k(\cdot,\Z) \rightarrow  C^k(\cdot,\R) \leftarrow \Om^k(\cdot).\]
Let $\eta^\Z\in C^3(\SU(2),\Z)$ be a smooth singular cocycle whose
image in $C^3(\SU(2),\R)$ is cohomologous to the image of $\eta$, and
let $\varpi_\pm^\Z \in C^2(\SU(2)_\pm,\Z)$ be primitives of the restriction
of $\eta^\Z$ to $\SU(2)_\pm$. Let $\sigma^\Z\in C^2(M,\Z)$ be such
that $\d\sigma^\Z=-k\Phi^*\eta^\Z$, and such that
$[(\sigma^\Z,k\eta^\Z)]\in H^3(\Phi,\Z)$ represents the lift of
$k[(\om,-\eta)]$ given by the pre-quantization. The upper map in the
commutative diagram above is given on $[(\sigma^\Z,k\eta^\Z)]$ by
\[ [(\sigma^\Z,k\eta^\Z)]\mapsto [\sigma^\Z_\pm+k\Phi^*\varpi^\Z_\pm].\]
Hence $[\sigma^\Z_\pm+k\Phi^*\varpi^\Z_\pm]\in H^2(M_\pm,\Z)$ are
integral lifts of $k[\om_{0,\pm}]$. Let $L_\pm\to M_\pm$ be
the corresponding $\SU(2)$-equivariant pre-quantum line bundles, so that 
\[ c_1(L_\pm)=[\sigma^\Z_\pm+k\Phi^*\varpi^\Z_\pm].\]

On the overlap, $M_\reg=M_+\cap M_-$, the difference between the
2-cocycles $\sigma^\Z_\pm+\Phi^*\varpi^\Z_\pm$ is
$k\Phi^*(\varpi^\Z_-|_{\SU(2)_\reg}-\varpi_\Z^+|_{\SU(2)_\reg})$. The
2-cochain $\varpi_-^\Z|_{\SU(2)_\reg}-\varpi_+^\Z|_{\SU(2)_\reg}\in
C^2(\SU(2)_\reg,\Z)$ is closed, and its cohomology class
is an integral lift of
$[\varpi_-|_{\SU(2)_\reg}-\varpi_+|_{\SU(2)_\reg}]=\Psi^*[\om_\O]\in H^2(\SU(2)_\reg,\R)$.  Hence it
represents the Chern class $c_1(J)$. We have shown that
\[c_1(L_-|_{M_\reg})-c_1(L_+|_{M_\reg})=k\Phi^* c_1(J)\]
and consequently $L_-|_{M_\reg}\cong L_+|_{M_\reg}\otimes \Phi^* J^k$.
Conversely, given a pair of pre-quantum line bundles $L_\pm$ with this
property, we may retrace the steps of this proof to obtain an integral
lift of $[k(\om,-\eta)]$.
\end{proof}

In particular, we see that if $(M,\om,\Phi)$ is pre-quantized at level
$k$, and $e$ is a regular value of $\Phi$, then the symplectic
quotient $M\qu \SU(2)$ inherits a level $k$ pre-quantization.  The
corresponding pre-quantum line bundle over $M\qu \SU(2)$ is $L_+\qu
\SU(2)=L_+|_{\Phinv}(e)/\SU(2)$ is a pre-quantum line bundle.

The pre-quantization result may be expressed in terms of Morita
trivializations:
\begin{proposition}\label{prop:preq}
  A level $k$ pre-quantization of a q-Hamiltonian $\SU(2)$-space
  $(M,\om,\Phi)$ gives rise to a Morita isomorphism,
\[ \C\simeq_\E \Phi^*\A^{k}.\]
\end{proposition}
\begin{proof}
  Pick Morita trivializations $ \C \simeq_{\F_\pm} \A^{k} $ over
  $\SU(2)_\pm$, with $\F_-\cong \F_+\otimes J^k$ on the overlap. The
  pre-quantum line bundles $L_\pm\to M_\pm$ defined by the level $k$
  pre-quantization satisfy $L_-\cong L_+\otimes\Phi^* J^k$ on the
  overlap.  Hence the Hilbert space bundles
\[ \E_\pm:=\on{Hom}_\C(L_\pm,\Phi^*\F_\pm)\]
(where $\on{Hom}_\C$ denotes continuous bundle homomorphisms) glue to
give the desired Morita isomorphism.
\end{proof}

\begin{proposition}
  Suppose $(M,\om,\Phi)$ is a level $k$ pre-quantized q-Hamiltonian
  $\SU(2)$-space. Assume that $e,c$ are regular values of $\Phi$.
  Then the imploded cross-section $(X,\om_X,\Phi_X)$ inherits a level
  $k$ pre-quantization.
\end{proposition}

\begin{proof}
Let $(M_\pm,\om_{0,\pm},\Phi_{0,\pm})$ carry the corresponding
pre-quantum line bundles with $L_-=L_+\otimes \Phi^* J^k$ on the
overlap. Since $X_\pm=(M_\pm\times\C^2)\qu \SU(2)$ are ordinary
Hamiltonian quotients, we obtain pre-quantizations of the Hamiltonian
$T$-spaces $(X_\pm,\om_X,\Phi_X)$. The pre-quantum line bundles 
$L_{X_\pm}$ satisfy $L_{X_\pm}|_Y\cong L_\pm|_Y$, hence
\[ L_{X_-}|_Y=L_{X_+}|_Y\otimes \Phi_Y^*J^k
=L_{X_+}|_Y\otimes \C_{k\rho}.
\]
We conclude that $L_{X_+}$ and $L_{X_-}\otimes \C_{-k\rho}$ patch to
define a global $T$-equivariant pre-quantum line bundle $L_X\to X$.
\end{proof}

\section{Quantization of q-Hamiltonian $\SU(2)$-spaces}
We are now in position to define the quantization of pre-quantized
q-Hamiltonian $\SU(2)$-spaces. We begin with a quick overview of the
quantization of ordinary \emph{Hamiltonian} $G$-spaces $(M,\om,\Phi)$.
Choose an invariant almost complex structure on $M$, compatible with
the symplectic form. Such an almost complex structure is unique up to
equivariant homotopy, and hence the isomorphism class of the resulting
equivariant $\Spin_c$-structure given by a $G$-equivariant spinor
bundle $\S$ is independent of this choice. We obtain a Morita
isomorphism $\Cl(TM)\simeq_{\S^{\on{op}}} \C$.  Given a pre-quantum line
  bundle $L\to M$, one can twist by $L$ to obtain a new
  $\Spin_c$-structure $\S\otimes L^{-1}$, hence a Morita
  isomorphism
\[ \Cl(TM)\simeq_{\S^{\on{op}}\otimes L}\C.\]
This allows us to define a push-forward map relative to 
$p\colon M\to \pt$,
\[ p_*\colon K_0^G(M,\Cl(TM))\to K_0^G(\pt)=R(G),\] 
and to set $\ca{Q}(M)=p_*([M])\in R(G)$. (For $G=\{e\}$, this is just
an integer.) Equivalently, $\ca{Q}(M)$ may be viewed as the
equivariant index of the $\Spin_c$-Dirac operator for the
$\Spin_c$-structure $\S\otimes L^{-1}$.  The quantization procedure for
Hamiltonian $G$-spaces is compatible with products:
\begin{equation}\label{eq:product}
 \ca{Q}(M_1\times M_2)= \ca{Q}(M_1)\ca{Q}(M_2).\end{equation}
For any $g\in G$, the value of the equivariant index $\ca{Q}(M)$ at
$g$ may be computed by Atiyah-Segal's localization theorem. On the
other hand, one has the Guillemin-Sternberg \emph{quantization
commutes with reduction} property: Let $\ca{Q}(M)^G\in \Z$ be the
multiplicity with which the trivial representation occurs in
$\ca{Q}(M)$. Then \cite{me:sym,me:si}
\[ \ca{Q}(M)^G=\ca{Q}(M\qu G).\]
Here the index $\ca{Q}(M\qu G)$ is well-defined if $0$ is a regular
value of $\Phi$ and the $G$-action on $\Phinv(0)$ is free. If the
action is only locally free, then $M\qu G$ is an orbifold and the
quantization is defined by the index theorem for orbifolds. In the
general case, if $0$ is not a regular value and $M\qu G$ is a singular
space, $\ca{Q}(M\qu G)$ may be defined by partial desingularization of
the singular symplectic quotient \cite{me:si}. 

Suppose now that $(M,\om,\Phi)$ is a compact q-Hamiltonian
$\SU(2)$-space, pre-quantized at level $k$. By combining the Morita
isomorphisms $\Phi^*\A^{2} \simeq_{\S} \Cl(TM)$
from Proposition \ref{prop:spinc} 
and $\C\simeq_\E \Phi^*\A^{k}$ from Proposition \ref{prop:preq} we
obtain a Morita isomorphism
\[ \Cl(TM)\simeq_{\S^{\on{op}}\otimes\E}\Phi^*\A^{k+2}.\]
This defines a push-forward map in $K$-homology, 
\[ K_0^{\SU(2)}(M,\Cl(TM))\to K_0^{\SU(2)}(\SU(2),\A^{k+2})\cong R_k(\SU(2)).\]
\begin{definition}
Let $(M,\om,\Phi)$ be a compact q-Hamiltonian $\SU(2)$-space,
pre-quantized at level $k$. We define the  \emph{quantization} 
$\ca{Q}(M)\in R_k(\SU(2))$ to be the push-forward of the $K$-homology fundamental class
$[M]\in K_0^{\SU(2)}(M,\Cl(TM))$, 
\[ \ca{Q}(M)=\Phi_*([M]).\]
\end{definition}
The properties of this quantization procedure for q-Hamiltonian spaces
are very similar to that for the Hamiltonian case: In particular, the
analogue to the `quantization commutes with products' property
\eqref{eq:product} holds, with the left hand side involving the fusion
product of q-Hamiltonian spaces, and the right hand side the product
in $R_k(\SU(2))$. However, while \eqref{eq:product} is rather obvious
in the Hamiltonian theory, its q-Hamiltonian counterpart is a
non-trivial fact (proved in \cite{al:for}).  In what follows, we will
focus on `localization' and `quantization commutes with reduction' for
q-Hamiltonian $\SU(2)$-spaces.

\section{Localization}
We had mentioned in \ref{subsec:two} that any $\tau\in R_k(\SU(2))$ is
determined by its values $\tau(t)$ at elements $t\in
T_{k+2}^\reg$. For a level $k$ pre-quantized q-Hamiltonian
$\SU(2)$-space $(M,\om,\Phi)$, the number $\ca{Q}(M)(t)$ may be
computed by localization to the fixed point set $M^t$ of $t$. By
equivariance, and since $t$ is regular, the moment map takes the fixed point set to the maximal
torus $T=\SU(2)^t$.
\begin{proposition}
  The restriction $\A^{k+2}|_T$ admits a $T_{k+2}$-equivariant Morita
  trivialization, 
\[ \C\simeq_\G \A^{k+2}|_T.
\]
This Morita trivialization is uniquely determined (up to equivalence)
by requiring that $\G|_e$ extends to an $\SU(2)$-equivariant Morita
trivialization of $\A^{k+2}|_e$.
\end{proposition}
\begin{proof}
  Choose $\SU(2)$-equivariant Morita trivializations $\C
  \simeq_{\F_\pm} \A^{k+2}|_{\SU(2)_\pm}$ such that on the overlap,
  $\F_-\cong \F_+\otimes J^{k+2}$. Restrict to $T$-equivariant Morita
  trivializations over
\[ T\cap \SU(2)_+ =T_{(-\rho,\rho)},\ \ T\cap \SU(2)_-=T_{(0,2\rho)}.\]
The intersection $T_{(-\rho,\rho)}\cap T_{(0,2\rho)}$ has two connected
components, $T_{(0,\rho)}$ and $T_{(\rho,2\rho)}$. 
The restrictions of $J^{k+2}$ to the two components are
\[ \begin{split}
J^{k+2}|_{T_{(0,\rho)}}&=T_{(0,\rho)}\times \C_{(k+2)\rho},\\   
J^{k+2}|_{T_{(\rho,2\rho)}}&=T_{(\rho,2\rho)}\times \C_{-(k+2)\rho}.\end{split}\]
Let 
\[ \G_+=\F_+|_{T_{(-\rho,\rho)}},\ \ \G_-=\F_-|_{T_{0,2\rho}}\otimes
\C_{(k+2)\rho}.\]
Then $\G_-\cong \G_+$ over $T_{(0,\rho)}$, while $\G_-=\G_+\otimes
\C_{2(k+2)\rho}$ over $T_{(\rho,2\rho)}$. 
But $T_{k+2}$ is exactly the subgroup of $T$ acting trivially on
$\C_{2(k+2)\rho}$.  That is, the bundles
$\G_\pm$ glue to define a $T_{k+2}$-equivariant Morita trivialization
\[ \C\simeq_\G \A^{k+2}|_T.\] 
By construction, $\G|_e$ extends to the unique (up to equivalence)
$\SU(2)$-equivariant trivialization $\F_+|_e$ of $\A|_e$.  Any other
$T_{k+2}$-equivariant Morita trivialization differs from $\G$ by twist
with a $T_{k+2}$-equivariant line bundle. Since $\dim T=1$ we have
$H^2_{T_{k+2}}(T)=H^2_{T_{k+2}}(\pt)$, hence such a line bundle is
detected by its restriction to $e$. Since only the trivial
$T_{k+2}$-representation extends to an $\SU(2)$-representation, the
proof is complete.
\end{proof}

\begin{remark}
The last part of the proof relied on $\dim T=1$. Indeed, the
corresponding statement for higher rank groups is more tricky
\cite{al:for}.
\end{remark}

\begin{proposition}\label{prop:spp}
  Suppose $\Phi\colon M\to \SU(2)$ is an equivariant map, and that we
  are given an equivariant Morita isomorphism
  $\Cl(TM)\simeq_\E\Phi^*\A^{k+2}$.  Then, for all regular elements
  $t\in T\cap \SU(2)_\reg$, and any component of the fixed point set
  $F\subset M^t$, the restriction $TM|_F$ inherits a distinguished
  $T_{k+2}$-equivariant $\Spin_c$-structure.
\end{proposition}
\begin{proof}
  By equivariance, and since $t$ is regular, $\Phi$ restricts to a map
  $\Phi_F\colon F\to \SU(2)^t=T$. Hence we have $T_{k+2}$-equivariant Morita
  isomorphisms 
\[ \C\simeq_{\Phi^*\G}
\Phi^*(\A^{k+2}|_T)\simeq_{\E^{\on{op}}|_F}\Cl(TM|_F).\] 
But a Morita trivialization of a Clifford algebra bundle is equivalent
to a $\Spin_c$-structure.
\end{proof}
Let $\L_F\to F$ be the $\Spin_c$-line bundle associated to this
$\Spin_c$-structure on $TM|_F$. 
\begin{remark}
  The line bundle $\L_F$ may be described as follows.  From
  $\Cl(TM)\simeq_\E\Phi^*\A^{k+2}$ we obtain a Morita trivialization,
\[ \C\simeq \Cl(TM)\otimes\Cl(TM)\simeq_{\E\otimes\E} \Phi^*\A^{2k+4}.\]
Over $M_\pm$, we have another Morita trvialization of
$\Phi^*\A^{2k+4}$ coming from the defining Morita trivializations of
$\A$ over $U_\pm$. The two Morita trivializations are related by line
bundles $\L_\pm\to M_\pm$, with $\L_-=\L_+\otimes \Phi^*J^{-(2k+4)}$ on
the overlap. The restriction of $J^{2k+4}$ to $T$ is
$T_{k+2}$-equivariantly trivial, and $\L_F$ is the
$T_{k+2}$-equivariant line bundle obtained by gluing $\L_\pm|_{F\cap
  M_\pm}$.
\end{remark}

Using Proposition \ref{prop:spp} we see that even though $M$ does not
come with a $\Spin_c$-structure, the fixed point contributions from
the usual Atiyah-Segal-Singer theorem \cite{at:1,at:2,at:3} are
well-defined. Indeed one has,
\begin{theorem}[Localization]\label{th:loc}
  Suppose $(M,\om,\Phi)$ is a compact q-Hamiltonian $\SU(2)$-space,
  pre-quantized at level $k$. For all $t\in T_{k+2}^\reg$, the number
  $\ca{Q}(M)(t)$ is given as a sum of fixed point contributions,
\[ \ca{Q}(M)(t)=\sum_{F\subset M^t}\ca{Q}(\nu_F)(t),\]
where $\ca{Q}(\nu_F)(t)$ is defined using the $T_{k+2}$-equivariant 
$\Spin_c$-structure on $TM|_F$. 
\end{theorem}

The proof of Theorem \ref{th:loc} is parallel to the proof of the
localization formula in Atiyah-Segal \cite{at:2}; details will be
given in \cite{al:for}. In the cohomological form of the index theorem,
the fixed point contributions $\ca{Q}(\nu_F)$ are given as integrals of certain
characteristic classes over $F$ (cf. \cite{du:he,al:fi}) 
\[ \ca{Q}(\nu_F)(t)=(\sigma(\L_F)(t))^{1/2} \int_F
\f{\wh{A}(F)\exp(\hh c_1(\L_F)) }{D_\R(\nu_F,t)}.\]
Here $\wh{A}(F)$ is the $\wh{A}$-class, and $D_\R(\nu_F,t)$ 
is given on the level of differential forms by  
\[ D_\R(\nu_F,t)=e^{\f{i\pi}{4}\on{rank}_\R(\nu_F)} {\det}_\R^{1/2}(1-t^{-1}e^{\f{1}{2\pi}\on{curv}_\R(\nu_F)}),\]
with $\on{curv}_\R(\nu_F)\in \Om^2(F,\mf{o}(\nu_F))$ the curvature
form for an invariant Riemannian connection. The expression in
parentheses lies in $\Om(F,\on{End}(\nu_F))$, with zeroth order term
the identity, and the (positive) square root of its determinant is
well-defined.  Finally $\L_F$ is the line bundle associated to the
$\Spin_c$-structure on $TM|_F$, the phase factor $\sigma(\L_\F)(t)\in
\U(1)$ is given by the action of $t$ on $\L|_F$, and
$\sigma(\L_\F)(t)^{1/2}$ is a suitable choice of square root.
\footnote{The square root is determined as follows. Let $\ca{S}_x$ be
  the fiber of the spinor module at any given $x\in F$.  Choose a
  $T_{k+2}$-invariant complex structure on $T_xM$, compatible with the
  orientation. Let $c_1,\ldots,c_{n/2}\in\U(1)$ be the eigenvalues
  (with multiplicities) for the action of $t$ on $T_xM$, and $u\in
  \U(1)$ the action of $t$ on the line
  $\on{Hom}_{\Cl(T_xM)}(\wedge_{\C} T_xM,\ca{S}_x)$. Then
\[ \sigma(\L_F)(t)^{1/2}=u\prod_{c_r\not=1} c_r^{1/2},\]
using the square roots of $c_r\not=1$ with positive imaginary part.}
If $F\subset M_+$, the $\Spin_c$-structure on $TM|_F$ is defined by
the almost complex structure on $M_+$, twisted by the line bundle
$L_+$.  Hence, the fixed point contribution can be written in
`Riemann-Roch' form:
\[ \ca{Q}(\nu_F)(t)=\sigma(L_+|_F)(t)\int_F\f{\on{Td}(F) \on{ch}(L_+|_F)}{D(\nu_{F,+},t)},\]
where $D(\nu_{F,+},t)$ is the equivariant characteristic class
\[ D(\nu_{F,+},t)=\on{det}_\C(1-t^{-1} e^{\f{i}{2\pi}\on{curv}_\C(\nu_{F,+})}),\]
with $\on{curv}_\C(\nu_{F,+})$ the curvature form for an invariant
Hermitian connection, and $\sigma(L_+|_F)(t)$ the phase factor defined
by the action of $t$ on $L_+|_F$. There is a similar formula for the case
$F\subset M_-$: 
\[ \ca{Q}(\nu_F)(t)=-t^{(k+2)\rho}\sigma(L_-|_F)(t)\int_F\f{\on{Td}(F) \on{ch}(L_-|_F)}{D(\nu_{F,-},t)}.\]
If $t=j(q^s)$ with $s=1,\ldots,k+1$, we have 
\[ -t^{(k+2)\rho}=(-1)^{s-1}.\]
This sign factor may be traced back to our choice of Morita
trivialization of $\A^{k+2}|_T$, which was chosen to be compatible
with the $\SU(2)$-equivariant Morita trivialization of $\A^{k+2}|_e$
(rather than that of $\A^{k+2}|_c$).

\begin{remark}
  A detailed check of the equivalence of the `$\Spin_c$' and 
  `Riemann-Roch' forms of the fixed point contribution may be found in
  \cite[Section 2.3]{al:fi}. In general, it is quite possible that $F$
  is contained neither in $M_+$ nor in $M_-$: this happens for
  instance for $M=D(\on{SO}(3))$, as discussed in the final Section of
  this paper.
\end{remark}

\begin{remark}
  The right hand side of the localization formula appears in
  \cite{al:fi}, as a `working definition' of the quantization of a
  q-Hamiltonian space.  However, in \cite{al:fi} it was not understood
  how to view this expression as the localization of an appropriate
  equivariant object on $M$.
\end{remark}

\section{Quantization commutes with reduction}
Suppose $(M,\om,\Phi)$ is a compact q-Hamiltonian $\SU(2)$-space, with
a pre-quantization at level $k$. For each $l=0,\ldots,k$, let $\Co_l$
be the conjugacy class of the element $\exp(\f{l}{k}\rho)$. If $\SU(2)$
acts freely (resp. locally freely) on $\Phinv(\Co_l)$, then 
\[ M\qu_{\Co_l}\SU(2)=(M\times \Co_l)\qu \SU(2)\cong \Phinv(\Co_l)/\SU(2)\] 
is a smooth
symplectic manifold (resp. orbifold), with a level $k$ pre-quantization from $M$. The
Riemann-Roch numbers 
\[ \ca{Q}(M\qu_{\Co_l}\SU(2))\in \Z\]
are thus defined. If $\SU(2)$ does not act locally freely, it is still
possible to define the Riemann-Roch numbers using a partial
desingularization, as in \cite{me:si}. 
\begin{theorem}[q-Hamiltonian quantization commutes with reduction]
\label{th:quantred}
  Let $(M,\om,\Phi)$ be a level $k$ pre-quantized q-Hamiltonian
  $\SU(2)$-manifold, and $\ca{Q}(M)\in R_k(\SU(2))$ its quantization.
  Let $N(l)\in\Z$ be the multiplicity of $\tau_l$ in $\ca{Q}(M)$. Then
\[ N(l)=\ca{Q}(M\qu_{\Co_l}\SU(2))\]
where the right hand side denotes the level $k$ quantization of the
symplectic quotient. 
\end{theorem}
A general proof of this result, for arbitrary simply connected groups,
can be found in \cite{al:fi}. Here we will present a much simpler
approach for the rank $1$ case. It is modeled after a similar proof
for the Hamiltonian case \cite[Appendix]{me:sym}.
\begin{proposition}
  Let $(M,\om,\Phi)$ be a level $k$ pre-quantized q-Hamiltonian
  $\SU(2)$-space. Suppose $\SU(2)$ acts (locally) freely on $\Phinv(e), \Phinv(c)$, so that
  the imploded cross-section $(X,\om_X,\Phi_X)$ is a smooth 
  Hamiltonian $T$-space, with a pre-quantization at level $k$. 
  Let $N_X(l),\,l\in \Z$ be the multiplicity function for the
  Hamiltonian $T$-space $X$, and $N(l),\ 0\le l\le k$ that for the
  q-Hamiltonian $\SU(2)$-space $M$. Then 
\[ N_X(l)=\begin{cases} N(l)&\mbox{ if }0\le l\le k\\
0&\mbox{ otherwise}\end{cases}\]
\end{proposition}
\begin{proof}
  We will only consider the case that $\SU(2)$ acts freely on
  $\Phinv(e),\Phinv(c)$.  The fact that $N_X(l)$ vanishes unless $0\le
  l\le k$ is an easy special case of the Hamiltonian `quantization commutes with
  reduction' theorem -- see e.g. \cite{du:sy}. The statement is thus
  equivalent to showing that $\ca{Q}(M)$ is the image, under the
  induction map $R_k(T)\to R_k(\SU(2))$, of $t^\rho \ca{Q}(M)(t)\in
  R(T)$ (restricted to $T_{k+2}$). That is, we have to show that for
  all $t=j(z)$, with $z\in \{q,q^2,\ldots,q^{k+1}\}$,
\[ 
\ca{Q}(M)(t)=\f{t^\rho\,\ca{Q}(X)(t)-t^{-\rho}\,\ca{Q}(X)(t^{-1})}{t^\rho-t^{-\rho}}
=\f{\ca{Q}(X)(t)}{1-t^{-2\rho}}+\f{\ca{Q}(X)(t^{-1})}{1-t^{2\rho}}.\]
The equivariant index theorem expresses $\ca{Q}(M)(t)$ as a sum of
fixed point contributions, $\ca{Q}(\nu_F)(t)$, as explained above.
Since $\SU(2)$ acts freely on $\Phinv(e),\Phinv(c)$, the fixed point
manifolds $F$ are all contained in $M_\reg$, hence we may work with
the Riemann-Roch form of teh fixed point contributions. By regularity,
$\Phi(F)\subset T^\reg$. Thus, either $F\subset Y$, or the image of
$F$ under the Weyl group action lies in $Y$. That is, all fixed point
manifolds come in pairs $F,F'$, with $F\in Y$ and $F'$ its image under
the action of the non-trivial Weyl group element. We have,
\[\ca{Q}(\nu_{F'})(t)=\ca{Q}(\nu_F)(t^{-1}).\]
Now, since $F\subset Y$ it also appears as a fixed point set in $X$.
The normal bundle of $F$ in $M$ splits as a direct sum of its normal
bundle $\nu_F^X$ in $X$ and the normal bundle of $Y$ in $M$, the
latter being $T$-equivariantly isomorphic to $\C_\alpha=\C_{2\rho}$.
Hence, the fixed point contributions are related by 
\[ \ca{Q}(\nu_F)(t)=\f{\ca{Q}(\nu_F^X)(t)}{1-t^{-2\rho}},\ \ 
\ca{Q}(\nu_{F'})(t)=\f{\ca{Q}(\nu_F^X)(t^{-1})}{1-t^{2\rho}}
.\]
Summing over all fixed point components $F\subset Y^t$, one obtains
all contributions to the fixed point formula for $X$, \emph{except} the
contributions from $F=M\qu \SU(2)$ and $F=M\qu_c \SU(2)$. From the explicit description of the normal
bundle of $M\qu \SU(2)$ as $\Phinv(0)\times_{\SU(2)}\C^2$, and the
identity, for $\xi\in\mf{su}(2)$,
\[ \det(1-z^{-1}e^{-\xi})=z^{-2}\det(1-z e^\xi)=z^{-2}\det(1-z
e^{-\xi})\]
we obtain, 
\[ D(\nu_{M\qu \SU(2)}^X,z^{-1})=z^{-2} D(\nu_{M\qu \SU(2)}^X,z).\]
Hence, the two terms for $F=M\qu \SU(2)$ cancel in the fixed point
formula for $X$. Similarly, the two contributions from $F=M\qu_c
\SU(2)$ cancel.
\end{proof}

\begin{proof}[Proof of Theorem \ref{th:quantred}]
  We have seen that $N(l)=N_X(l)$.  From the `quantization commutes
  with reduction theorem' for Hamiltonian $\U(1)$-spaces \cite{du:sy},
  we know that $N_X(l)$ is the Riemann-Roch number of the level $k$
  quantization of a symplectic quotient of $X$:
\[ N_X(l)=\ca{Q}(\Phi_X^{-1}(\f{i\pi l}{k})/\U(1))
=\ca{Q}(M\qu_{\Co_l}\SU(2)). \qedhere\]
\end{proof}

One obtains the multiplicities $N(l)$ by the orthogonality relations
\eqref{eq:orth}.
Writing $N(l)=\ca{Q}(M\qu_{\Co_l} \SU(2))$ we obtain,
\[ \ca{Q}(M\qu_{\Co_l}\SU(2))=\sum_{s=1}^{k+1} \f{|q^s-q^{-s}|^2}{2k+4}\,\,\tau_l(j(q^s))\,\,\ca{Q}(M)(j(q^s))
.\]

\section{Examples}
Using the localization formula, we can compute the quantizations
$\ca{Q}(M)\in R_k(\SU(2))$ for our basic examples. 
Recall that $\tau_n,\ n=0,\ldots,k$ are the
basis elements of $R_k(\SU(2))$.

\subsection{The double}\label{subsec:double}
We begin with the q-Hamiltonian $\SU(2)$-space $D(\SU(2))$.  Recall
that this space is pre-quantizable at any integer level $k\ge
1$. 
\begin{proposition}\label{prop:dsu2}
The level $k$ quantization of the double $D(\SU(2))$ is given by 
\[ \ca{Q}(D(\SU(2))=\sum_{j=0}^{[\f{k}{2}]} (k+1-2j) \tau_{2j}.\]
Here $[x]$ denotes the largest integer less than or equal to
$x$. Equivalently, 
\[ \ca{Q}(D(\SU(2))(j(q^s))=\f{2k+4}{ |q^s-q^{-s}|^{2}}\]
for $s=1,\ldots,k+1$.
\end{proposition}
\begin{proof}
We first verify the equivalence of the two formulas. Using the known
formulas for products of $\tau_n$'s, one finds
that
\[ \sum_{j=0}^{[\f{k}{2}]} (k+1-2j) \tau_{2j}=\sum_{n=0}^k (\tau_n)^2.\] 
Write $z=q^s$. Then
\[\begin{split} 
\sum_{n=0}^k (\tau_n(j(z)))^2&=-\f{1}{|z-z^{-1}|^2}
\sum_{n=0}^k (z^{n+1}-z^{-(n+1)})^2\\&= 
\f{1}{|z-z^{-1}|^2}\sum_{n=0}^k (2-z^{2(n+1)}-z^{-2(n+1)})=
\f{2k+4}{|z-z^{-1}|^2},\end{split}\]
where the sum is evaluated as a geometric series (using
$z^{k+2}=(-1)^s$). We next compare this result to the fixed point
computation for $M=D(\SU(2))$ (the following computation may be found
in \cite{al:fi}). Since the action of $\SU(2)$ on
$M=\SU(2)\times\SU(2)$ is by conjugation on each factor, and $j(z)$ is
a regular element, its fixed point set is
\[ M^{j(z)}=T\times T=:F.\] 
Note that $\Phi(F)=\{e\}$, in particular $F\subset M_+$. The induced
symplectic structure on $F$ is the standard symplectic structure on
$T\times T$, defined by the inner product:
\[ \om_F=\pr_1^*\theta_T\cdot \pr_2^*\theta_T\]
where $\pr_i\colon T\times T\to T$ are the two projections. 
The symplectic volume of $F$ is 
\[ \vol(F)=\int_{T\times T}\om_F=(\int_T \theta_T)\cdot (\int_T \theta_T)=\alpha\cdot\alpha=2.\]
The $\Spin_c$-line bundle $\L_F$ comes from the level $k+2$ Morita
isomorphism $\Cl(TM)\simeq \Phi^*\A^{k+2}$, 
\[ \C\simeq  \Cl(TM)\otimes\Cl(TM)\simeq \Phi^*\A^{2k+4}\]
hence it is isomorphic to the $2k+4$-th power of the level $1$
pre-quantum line bundle over $F$. (We are using that $H^2(M,\Z)=0$.)
Hence $\hh c_1(\L_F)=(k+2)\om_F$. By considering the action at
$x=(e,e)\in F$, one checks that $\zeta(\L_F)(t)=1$.  Indeed, the
$\Spin_c$-structure on $T_xM$ extends to an $\SU(2)$-equivariant
$\Spin_c$-structure, and the corresponding representation of $\SU(2)$
on $\L_F|_x$ is necessarily trivial.  The normal bundle to $F$ in $M$
is a trivial bundle 
\[ \nu_F=\su(2)/\t\oplus \su(2)/\t=\C\oplus \C^-,\] 
with $T$ acting by
weight $2$ on the first summand and $-2$ on the second summand. Hence
\[ 
\f{\zeta_F(t)^{1/2}}{D_\R(\nu_F,t)}=\f{1}{|(1-z^{2})(1-z^{-2})|}=\f{1}{|z-z^{-1}|^2}.\]
Since finally $\wh{A}(F)=1$, the fixed point contribution is 
\[ \chi(\nu_F,j(z))=\int_F \f{e^{(k+2)\om_F}}{|z-z^{-1}|^2}
=\f{2k+4}{|z-z^{-1}|^2},\]
as claimed.
\end{proof}
Recall now that $M(\Sig_h)=D(\SU(2))^h\qu\SU(2)$ is the moduli space
of flat $\SU(2)$-bundles over a surface of genus $h$. Using that
quantization commutes with products, we have
$\ca{Q}(D(\SU(2))^h)=\ca{Q}(D(\SU(2)))^h$.  Together with the
quantization commutes with reduction principle we hence obtain the
Verlinde formula for this moduli space (cf.  \cite{sz:ver}):
\[ \ca{Q}(M(\Sig_h))=\sum_{s=1}^{k+1}
\Big(\f{|q^s-q^{-s}|^2}{2k+4}\Big)^{1-h}
=\sum_{s=1}^{k+1}
\Big(\f{2\sin^2(\f{s\pi}{k+2})}{k+2}\Big)^{1-h}
.\]

\subsection{Conjugacy classes}
We had seen that the conjugacy classes $\Co\subset \SU(2)$
admitting a level $k$ pre-quantizations are precisely those of
elements $\exp(\f{n}{k}\rho)$ with $0\le n\le k$.
\begin{proposition}
The level $k$ quantization of the conjugacy class $\Co=\SU(2).\exp(\f{n}{k}\rho)$
is given by 
\begin{equation}\label{eq:conj1} 
\ca{Q}(\Co)=\tau_n.\end{equation}
Equivalently, for $s=1,\ldots,k+1$, 
\begin{equation}\label{eq:conjugacy class}
 \ca{Q}(\Co)(j(q^s))=\f{q^{s(n+1)}-q^{-s(n+1)}}{q^s-q^{-s}}.
\end{equation} 
\end{proposition}
\begin{proof}
  The equivalence of the two formulations follows from the discussion
  in Section \ref{subsec:two}.  Write $z=q^s$. If $n<k$, then
  $\Phi(\Co)\subset \SU(2)_+$. The symplectic form on $\Co=\Co_+$
  identifies $\Co$ with the coadjoint orbit of $\f{n}{k}\rho$, and the
  level $k$ pre-quantization corresponds to the usual (level 1)
  pre-quantization of the orbit through $n\rho$.  Written in
  Riemann-Roch form, the fixed point contributions for the conjugacy
  class are just the same as those for the coadjoint orbit, given by
  \eqref{eq:conjugacy class}. If $n=k$, the conjugacy class $\Co$
  coincides with the central element $\{c\}$.  Since $z^{k+2}=(-1)^s$
  we have,
\[ \chi_k(z)=\f{z^{k+1}-z^{-(k+1)}}{z-z^{-1}}
            =\f{z^{k+2}-z^{-(k+2)}z^2}{z^2-1}=-(-1)^s\]
which on the other hand is also the fixed point contribution for
$\ca{Q}(\Co)(j(z))$, for $\Co\in \Phinv(\SU(2)_-)$. This
gives \eqref{eq:conjugacy class} for $n=k$. 
\end{proof}

As a consequence, we may compute the level k quantization of 
\[ M(\Sig_h^r;\Co_1,\ldots,\Co_r)=D(\SU(2))^h\times\Co_1\times\cdots\times\Co_r\]
where $\Co_i,\ i=1,\ldots,r$ are conjugacy classes of
elements $\exp(\f{l_i}{k}\rho)$ with $0\le l_i\le k$. One obtains, 
\[ 
\ca{Q}(M(\Sig_h^r;\Co_1,\ldots,\Co_r))=\sum_{s=1}^{k+1}
\Big(\f{|q^s-q^{-s}|^2}{2k+4}\Big)^{1-h} \tau_{l_1}(q^s)\cdots
\tau_{l_r}(q^s).\] For $h=0$ and $r=3$, the right hand side of this
formula are the fusion coefficients. That is, 
\[ \ca{Q}(M(\Sig_0^3:\Co_1,\Co_2,\Co_3))=N^{(k)}_{l_1,l_2,l_3}.\]

\subsection{The 4-sphere}
Recall that the q-Hamiltonian space $S^4$ admits a unique
pre-quantization for all $k$. 
\begin{proposition}
The level $k$ quantization of the 4-sphere is given by 
\[ \ca{Q}(S^4)=\sum_{n=0}^k \tau_n.\]
Equivalently, for $s=1,\ldots,k+1$
\[ \ca{Q}(S^4)(j(q^s))=\begin{cases}2\ |1-q^{-s}|^{-2},\ \ & s \mbox{ odd}\\
0 & s \mbox{ even}\\
\end{cases}\]
\end{proposition}
\begin{proof}
Write $z=q^s$. We first verify the equivalence of the two formulas:  
\[ \begin{split}
\sum_{n=0}^k \tau_n(j(q(z))
&=\f{1}{z-z^{-1}}\sum_{n=0}^k(z^{n+1}-z^{-(n+1)})
\\&=\f{1}{z-z^{-1}}(\f{z-z^{k+2}}{1-z}-\f{z^{-1}-z^{-(k+2)}}{1-z^{-1}}).
\end{split}\]
If $s$ is even, then $z^{k+2}=1$ and the two terms cancel. If $s$ is
odd, then $z^{k+2}=-1$ and we obtain, writing $(z-z^{-1})=(1-z^{-1})(z+1)$,  
that 
\[ \sum_{n=0}^k \tau_n(j(z))=\f{2} {(1-z^{-1})(1-z)}=\f{2}{|1-z^{-1}|^2}.\]  
The fixed point set of $t$ consists
  of the `north pole' $\Phinv(e)$ and the `south pole' $\Phinv(c)$. By
  construction, $S^4_\pm$ are identified with open balls in $\C^2$,
  with the standard $\SU(2)$-action. Hence the weights for the
  $T\subset \SU(2)$-action are $+1,-1$ respectively, and the fixed
  point formulas give (using $j(z)^{(k+2)\rho}=z^{k+2}=(-1)^s$)
\[
\ca{Q}(S^4)(j(z))=\f{1}{(1-z)(1-z^{-1})}-(-1)^s\f{1}{(1-z)(1-z^{-1})},\]
as needed. 
\end{proof}

\subsection{Moduli spaces of flat $\SO(3)$-bundles}
The symplectic quotient
\[ D(\SO(3))^h\qu \SO(3)\] 
of an $h$-fold product of $D(\SO(3))$'s (viewed as
$q$-Hamiltonian $\SO(3)$-spaces) is the moduli space of
flat $\SO(3)$-bundles over a surface of genus $h$. It has two 
connected components, given as symplectic quotients of 
$D(\SO(3))^h$ where $D(\SO(3))$ is now viewed as a q-Hamiltonian
$\SU(2)$-space: 
\begin{equation}\label{eq:so3} D(\SO(3))^h\qu
\SO(3)=D(\SO(3))^h\qu \SU(2)\ \cup\ D(\SO(3))^h\qu_c \SU(2).
\end{equation}
The two components correspond
to the trivial and the non-trivial $\SO(3)$-bundle over the surface.  To
obtain Verlinde numbers for these moduli spaces, we need to work out
the quantization of the q-Hamiltonian $\SU(2)$-space $D(\SO(3))$.

We had seen that $D(\SO(3))$ is pre-quantizable at level $k$ if and
only if $k$ is even. The different pre-quantizations are a principal
homogeneous space under the torsion subgroup of $H^2(D(\SO(3)),\Z)$.
In fact this group is all torsion, and
  \[\begin{split}
    H^2(D(\SO(3)),\Z)&=H^2_{\Z_2\times\Z_2}(D(\SU(2)),\Z)\\
                     &=H^2_{\Z_2\times\Z_2}(\pt,\Z)\\
                     &=\on{Hom}(\Z_2\times\Z_2,\U(1)) 
\end{split}.\]
Letting $\C_\phi$ denote the 1-dimensional representation given by
$\phi\in\on{Hom}(\Z_2\times\Z_2,\U(1))$, this group acts by tensoring
with the flat line bundle 
\[D(\SU(2))\times_{\Z_2\times\Z_2}\C_\phi.\]
Let $T'=T/\Z_2$ be the maximal torus in $\SO(3)$, and $N(T)\subset
\SU(2)$, $N(T')\subset \SO(3)$ the normalizers. Similarly, for
elements
$a,b,\ldots$ of $\SU(2)$ we denote by $a',b',\ldots$ their images in
$\SO(3)$. 
\begin{lemma}
For any $t\in T_{\on{reg}}\subset \SU(2)$, the fixed point set of its
action on $\SO(3)=\SU(2)/\Z_2$ is $T'=T/\Z_2$ unless $t^2=c$, in which
case it is $N(T')=N(T)/\Z_2$.   
\end{lemma}
\begin{proof}
  For $a\in \SU(2)$, the element $a'$ is fixed under $\Ad_t$ if and
  only if $a$ is fixed up to a central element, i.e.
  $tat^{-1}a^{-1}\in Z(\SU(2))$. If this central element is $e$, this
  just means $a\in T$. If the central element is $c$, then
  $at^{-1}a^{-1}=t^{-1}c$ shows that $a\in N(T)$ represents the
  non-trivial Weyl element $w$, and $c=tw(t^{-1})=t^2$. We have thus
  shown that the fixed point set of a regular element $t$ is the image
  of $T$ in $\SO(3)$, unless $t^2=c$ in which case it is the image of
  the normalizer $N(T)$.
\end{proof}

Let us consider the fixed contributions of any $t=j(q^s),\ 
s=1,2,\ldots,k+1$ for the q-Hamiltonian space $D(\SO(3))$, for $k$
even. Note that $t^2=c\Leftrightarrow s=k/2+1$, and so we have to 
consider two cases:   

Case 1: $s\not=1+\f{k}{2}$, i.e. $t^2\not=c$. Then
$D(\SO(3))^t=T'\times T'=:F$ is connected, and its moment map image is
$\{e\}$. Since $\SU(2)$ acts trivially on the fiber of $L_+$ at
$(e',e')\subset F$, the action of $t$ on $L_+|_F$ is trivial. Hence
the fixed point contribution is just $1/4$ that of the
corresponding fixed point manifold in $D(\SU(2))$:
\[ \chi(\nu_F,t)=\f{1}{4} \f{2k+4}{|q^s-q^{-s}|^2}
=\f{1}{4\sin^2(\f{\pi s}{k+2})}(\f{k}{2}+1)
.\]

Case 2: $s=1+\f{k}{2}$, i.e. $t^2=c$ and $q^s=i$. Then $D(\SO(3))^t
=N(T')\times N(T')$ has four connected components, indexed by the
elements of $u=(u_1,u_2)\in W\times W=\Z_2\times\Z_2$. Choose 
\[ n=\left(\begin{array}{cc}0&1\\-1&0\end{array}\right)\in N(T)\]
as a lift of the non-trivial Weyl group element, and let $n'\in N(T')$
its image. Then each fixed point component $F_u$ has a base point 
\[x_u\in \{(e',e'),(n',e'),(e',n'),(n',n')\}\] 
with the property $\Phi(x_u)=e$. For any given choice of the
pre-quantization, one finds that the contribution of the component
labeled by $u=(u_1,u_2)$ is of the form, \footnote{The computation is
  similar to that in Section \ref{subsec:double}. In particular, the
  symplectic volume of the 2-torus $F_u$ may be computed by working
  out $\om_{F_u}$ in coordinates; one finds $\vol(F_u)=1/2$. See
  \cite{al:ve} for more general calculations along these lines.}
\[  \chi(\nu_{F_u},t)=\f{\lambda(u)}{4} \f{2k+4}{|q^s-q^{-s}|^2}
=\f{\lambda(u)}{4} (\f{k}{2}+1) .\] 
where $\lambda(u)\in \U(1)$ is given by the action of $t$ on $L_+|_{m_u}$. For
$u=(1,1)$, this phase factor is $\lambda(u)=1$ as above. The total
fixed point contribution is obtained by summing over all
$u=(u_1,u_2)$:
\[ \ca{Q}(D(\SO(3))(q^{k/2+1})=(\f{k}{2}+1)\sum_u\f{\lambda(u)}{4}.\] 

Let $\chi\in R_k(\SU(2))$ be defined by 
\begin{equation}\label{eq:chi}
 \chi=\sum_{j=0}^{k/2} (-1)^j \tau_{2j}
=
\tau_0-\tau_2+\tau_4\cdots +(-1)^{k/2} \tau_k.\end{equation}
Using the orthogonality relations for level $k$ characters, one finds
that
\[ \chi(q^{k/2+1})=\f{k}{2}+1,\ \chi(q^s)=0\mbox{ for }s\not=k/2+1.\]
From the localization contributions, we see:  
\[ \ca{Q}(D(\SO(3)))=\f{1}{4}\Big(\ca{Q}(D(\SU(2)))
+\sum_{u\not=(1,1)} \lambda(u)\ \chi\Big).\]  
It remains to understand the sum $\sum_{u\not=(1,1)} \lambda(u)$.
\begin{lemma}
  For every even $k$, and any $\phi\in
  \on{Hom}(\Z_2\times\Z_2,\U(1))$, the space $D(\SO(3))$ admits a unique
  pre-quantization at level $k$ with the property that
  \[\lambda(u)=(-1)^{k/2}\phi(u)\]
  for all $u\not=(1,1)$.
\end{lemma}
\begin{proof}
Changing the pre-quantization by $\phi\in
\on{Hom}(\Z_2\times\Z_2,\U(1))$ changes $\lambda(u)$ to
$\ti{\lambda}(u)=\lambda(u)\phi(u)$. This shows uniqueness.
For existence, we have to find a pre-quantization with
$\lambda(u)=(-1)^{k/2}$ for $u\not=(1,1)$. In fact, it is 
enough to find such a pre-quantization for $k=2$. (The
  general case will then follow by taking the $k/2$-th power of the
  pre-quantization at level $2$.) 
  
  For $k=2$, and any of the four possible pre-quantizations, write 
\[ \ca{Q}(D(\SO(3)))=\sum_{l=0}^2 N(l)\tau_l.\]
The localization formulas for $q,q^2,q^3$ give equations
\[ \begin{split} N(0)+\sqrt{2} N(1)+N(2)&=1,\\
  N(0)-N(2)&=\hh +\hh \sum_{u\not=(1,1)}\lambda(u),\\
  N(0)-\sqrt{2} N(1)+N(2)&=1.
\end{split}\]
The first and third equation give $N(1)=0$ and $N(0)+N(2)=1$.  In
particular, $N(0)-N(2)$ is an odd integer. The second equation shows
that $\sum_{u\not=(1,1)} \lambda(u)$ is a real number.  A change of
pre-quantization produces a sign change of exactly two of the
$\lambda(u)$'s with $u\not=(1,1)$.  Since $\sum_{u\not=(1,1)}
\ti{\lambda}(u)$ is again a real number, it follows that all
$\lambda(u)$ are real, and hence equal to $\pm 1$. The number of
$\lambda(u)$'s equal to $-1$ must be odd, or else the second equation
would give that $N(0)+N(2)=0$ or $=2$, contradicting that $N(0)-N(2)$
is odd. Hence, either all three $\lambda(u)$'s with $u\not=(1,1)$ are
equal to $-1$, or exactly one of them $\lambda(u)$ equals $-1$ and the
other two are equal to $+1$.  The resulting four cases must correspond
to the four pre-quantizations. In particular, there is a unique level
$2$ pre-quantization such that $\lambda(u)=-1$ for all
$u\not=(-1,-1)$.
\end{proof}
Let $\delta_{\phi,1}$ be equal to $1$ if $\phi=1$, equal to $0$
otherwise.  Then $\sum_u \phi(u)=4\delta_{\phi,1}$, i.e.
$\sum_{u\not=(1,1)}\phi(u)=-1+4\delta_{\phi,1}$. It follows that
\[ \ca{Q}(D(\SO(3)))=\f{1}{4}\Big(\ca{Q}(D(\SU(2)))
+(-1)^{k/2}(-1+4\delta_{\phi,1})\ \chi\Big) .\] From the known
expansions of $\ca{Q}(D(\SU(2)))$ (Proposition \ref{prop:dsu2})
and $\chi$ (Equation \eqref{eq:chi}) in the basis $\tau_j$, we
finally obtain:
\begin{theorem}
For $k$ even, let $D(\SO(3))$ carry the level $k$ pre-quantization
labeled by $\phi\in\on{Hom}(\Z_2\times\Z_2,\U(1))$. Then 
\[ \ca{Q}(D(\SO(3)))=\f{1}{4}\sum_{j=0}^{k/2}\big(k+1-2j +(-1)^{j+k/2}(-1+4\delta_{\phi,1})\big)\tau_{2j}.\] 
Equivalently, for $s=1,\ldots,k+1$,
\[ \ca{Q}(D(\SO(3)))(j(q^s))=\begin{cases}\f{1}{4}\sin^{-2}(\f{\pi
      s}{k+2})(\f{k}{2}+1)&s\not=\f{k}{2}+1\\
\f{1}{4}\big(1+(-1)^{k/2}(-1+4\delta_{\phi,1})\big)(\f{k}{2}+1)
&s=\f{k}{2}+1
\end{cases}\]
\end{theorem}

\vskip.2in
Dividing into the various subcases, the formula reads,
\[
\ca{Q}(D(\SO(3)))=\begin{cases}(\f{k}{4}+1)\tau_0+(\f{k}{4}-1)\tau_2+\f{k}{4}\tau_4+(\f{k}{4}-2)\tau_6+\cdots
& \phi=1,\  k=0\mod 4\\
\f{k}{4}\tau_0+\f{k}{4}\tau_2+(\f{k}{4}-1)\tau_4+(\f{k}{4}-1)\tau_6+\cdots
&\phi\not=1,\  k=0\mod 4\\
(\f{k-2}{4})\tau_0+(\f{k-2}{4}+1)\tau_2+(\f{k-2}{4}-1)\tau_4+(\f{k-2}{4})\tau_6+\cdots
&\phi=1,\  k=2\mod 4\\
(\f{k-2}{4}+1)\tau_0+(\f{k-2}{4})\tau_2+(\f{k-2}{4})\tau_4+(\f{k-2}{4}-1)\tau_6+\cdots
&\phi\not=1,\  k=2\mod 4
\end{cases}
\]
Using this result, in combination with `quantization commutes with
reduction', it is now straightforward to compute the quantizations
(Verlinde numbers) for the moduli spaces \eqref{eq:so3}. Note that there
are many different pre-quantizations, since one can choose a different
$\phi$ for each factor. The case with boundary (markings) is still
more complicated, and will be discussed elsewhere.
\begin{remark}
  For $k=0\mod 4$, the result above was proved about eight years ago
  in joint work \cite{al:ve} with Anton Alekseev and Chris Woodward.
  Pantev \cite{pa:cm} and Beauville \cite{be:ve} had earlier obtained
  obtained similar results using techniques from algebraic geometry.
\end{remark}

\def\cprime{$'$} \def\polhk#1{\setbox0=\hbox{#1}{\ooalign{\hidewidth
  \lower1.5ex\hbox{`}\hidewidth\crcr\unhbox0}}} \def\cprime{$'$}
  \def\cprime{$'$} \def\polhk#1{\setbox0=\hbox{#1}{\ooalign{\hidewidth
  \lower1.5ex\hbox{`}\hidewidth\crcr\unhbox0}}} \def\cprime{$'$}
  \def\cprime{$'$}
\providecommand{\bysame}{\leavevmode\hbox to3em{\hrulefill}\thinspace}
\providecommand{\MR}{\relax\ifhmode\unskip\space\fi MR }
\providecommand{\MRhref}[2]{%
  \href{http://www.ams.org/mathscinet-getitem?mr=#1}{#2}
}
\providecommand{\href}[2]{#2}

\end{document}